\documentclass[aos,seceqn,nameyear,dvips]{arximspdf}
\usepackage{graphics}
\doi{10.1214/09-AOS789}
\volume{38}
\issue{4}
\pubyear{2010}
\firstpage{2388}
\lastpage{2421}

\makeatletter

\newtheorem{theorem}{Theorem}[section]
\newtheorem{proposition}{Proposition}[section]
\newtheorem{lemma}{Lemma}[section]

\makeatother

\begin{document}
\begin{frontmatter}

\title{Simultaneous nonparametric inference of time~series\thanksref{T1}}
\runtitle{Simultaneous nonparametric inference}

\thankstext{T1}{Supported in part by NSF Grant DMS-04-78704.}

\begin{aug}
\author[A]{\fnms{Weidong} \snm{Liu}\corref{}\ead[label=e1]{liuweidong99@gmail.com}} and
\author[B]{\fnms{Wei Biao} \snm{Wu}\ead[label=e2]{wbwu@galton.uchicago.edu}}
\runauthor{W. Liu and W. B. Wu}
\affiliation{University of Pennsylvania and University of Chicago}
\address[A]{Department of Statistics\\
University of Pennsylvania\\
3730 Walnut Street\\
Philadelphia, Pennsylvania 19104\\
USA\\
\printead{e1}}
\address[B]{Department of Statistics\\
University of Chicago\\
5734 S. University Avenue\\
Chicago, Illinois 60637\\
USA\\
\printead{e2}}
\end{aug}

\received{\smonth{8} \syear{2009}}
\revised{\smonth{12} \syear{2009}}

%
\begin{abstract}
We consider kernel estimation of marginal densities and regression
functions of stationary processes. It is shown that for a wide
class of time series, with proper centering and scaling, the
maximum deviations of kernel density and regression estimates are
asymptotically Gumbel. Our results substantially generalize
earlier ones which were obtained under independence or beta mixing
assumptions. The asymptotic results can be applied to assess
patterns of marginal densities or regression functions via the
construction of simultaneous confidence bands for which one can
perform goodness-of-fit tests. As an application, we construct
simultaneous confidence bands for drift and volatility functions
in a dynamic short-term rate model for the U.S. Treasury yield
curve rates data.
\end{abstract}

%

\end{frontmatter}

\section{Introduction}

Consider the nonparametric time series
regression model
%
%
\begin{equation}\label{eq:nreg}
Y_i = \mu(X_i)\,dt + \sigma(X_{i})\eta_i,
\end{equation}
where $\mu(\cdot)$ [resp., $\sigma^2(\cdot)$] is an unknown regression
(resp., conditional variance) function to be estimated, $(X_i,
Y_i)$ is a stationary process and $\eta_i$ are unobserved
independent and identically distributed (i.i.d.) errors with $\mathsf{E}
\eta_i = 0$ and $\mathsf{E}\eta^2_i = 1$. Let the regressor $X_i$ be a
stationarity causal process
%
%
\begin{equation}\label{a1}
X_i = G(\ldots, \varepsilon_{i-1},\varepsilon_i),
\end{equation}
where $\varepsilon_{i}$ are i.i.d. and the function $G$ is such that
$X_i$ exists. Assume that $\eta_i$ is independent of $(\ldots,
\varepsilon_{i-1},\varepsilon_i)$. Hence, $\eta_i$ and $(\mu(X_i),
\sigma(X_i))$ are independent. As a special case of
(\ref{eq:nreg}), a particularly interesting example is the
nonlinear autoregressive model
%
%
\begin{equation}\label{eq:nlts19}
Y_i = \mu(Y_{i-1}) + \sigma(Y_{i-1}) \eta_i,
\end{equation}
where $X_i = Y_{i-1}$ and $\varepsilon_i = \eta_{i-1}$. Many
nonlinear time series models are of form (\ref{eq:nlts19}) with
different choices of $\mu(\cdot)$ and $\sigma(\cdot)$. If the form
of $\mu(\cdot)$ is not known, we can use the Nadaraya--Watson
estimator
%
%
\begin{equation}\label{eq:nwe}
\mu_n(x)=\frac{1}{n b f_{n}(x)}
\sum_{k=1}^{n}K \biggl(\frac{X_{k}-x}{b} \biggr)Y_{k},
\end{equation}
where $K$ is a kernel function with $K(\cdot) \ge0$ and $\int
_{\mathsf{R}}
K(u) \,d u = 1$, the bandwidths $b = b_n \to0$ and $n b_n \to
\infty$, and
\[
f_{n}(x)=\frac{1}{nb}
\sum_{k=1}^{n} K \biggl(\frac{X_{k}-x}{b} \biggr)
\]
is the kernel density estimate of $f$, the marginal density of
$X_i$. Asymptotic properties of nonparametric estimates for time
series have been widely discussed under various strong mixing
conditions; see \citet{robin83}, \citet{gyorfietal89},
\citet{tjos94}, \citet{bosq96}, \citet{DL99} and
\citet{fan03}, among others.

Under appropriate dependence conditions [see, e.g.,
\citet{robin83}, \citet{WM02}, \citet{fan03} and
\citet{wu05}], we have the central limit theorem
\[
\sqrt{n b} [f_{n}(x)-\mathsf{E}f_{n}(x)] \Rightarrow
N(0, {\lambda_K f(x)})\qquad\mbox{where }
\lambda_K = \int_{\mathsf{R}}K^2(u) \,d u.
\]
The above result can be used to construct point-wise confidence
intervals of $f(x)$ at a fixed $x$. To assess shapes of density
functions so that one can perform goodness-of-fit tests, however,
one needs to construct \textit{uniform} or \textit{simultaneous
confidence bands} (SCB). To this end, we need to deal with the
maximum absolute deviation over some interval $[l, u]$:
%
%
\begin{equation}\label{eq:J91}
\Delta_n := \sup_{l\leq x\leq u}\frac{\sqrt{n b}}
{\sqrt{\lambda_K f(x)}}|f_{n}(x)-\mathsf{E}f_{n}(x)|.
\end{equation}
In an influential paper, \citet{BR73} obtained an
asymptotic distributional theory for $\Delta_n$ under the assumption
that $X_i$
are i.i.d. It is a very challenging problem to generalize their result to
stationary processes where dependence is the rule rather than the
exception. In their paper Bickel and Rosenblatt applied the very deep
embedding theorem of approximating empirical processes of independent random
variables by Brownian bridges with a reasonably sharp rate
[\citet{brill69}, Koml\'{o}s, Major and Tusn\'{a}dy (\citeyear
{KMT75}, \citeyear{KMT76})].
For stationary processes, however,
such an approximation with similar rates can be extremely difficult to
obtain. \citet{DP87} obtained a weak invariance principle for
empirical distribution functions. In 1998, Neumann (\citeyear{neum98}) made a breakthrough
and proved a very useful result for $\beta$-mixing processes whose mixing
rates decay exponentially quickly. Such processes are very weakly
dependent. For mildly weakly dependent processes, the asymptotic
problem of
$\Delta_n$ remains open. Fan and Yao [(\citeyear{fan03}), page 208]
conjectured that
similar results hold for stationary processes under certain mixing
conditions. Here we shall solve this open problem and establish an
asymptotic theory for both short- and long-range dependent processes.
It is
shown that, for a wide class of short-range dependent processes, we can
have a similar asymptotic distributional theory as \citet{BR73}.
However, for long-range dependent processes, the asymptotic
behavior can be sharply different. One observes the dichotomy
phenomenon: the
asymptotic properties depend on the interplay between the strength of
dependence and the size of bandwidths. For small bandwidths, the limiting
distribution is the same as the one under independence. If the
bandwidths are large, then the limiting distribution is half-normal [cf.
(\ref{eq:aug241})].

A closely related problem is to study the asymptotic uniform
distributional theory for the Nadaraya--Watson estimator
$\mu_n(x)$. Namely, one needs to find the asymptotic distribution
for $\sup_{x \in T} |\mu_n(x) - \mu(x)|$, where $T = [l, u]$. With
the latter result, one can construct an asymptotic $(1-\alpha)$
SCB, $0 < \alpha< 1$, by finding two functions $\mu_n^{\mathrm{lower}}(x)$
and $\mu_n^{\mathrm{upper}}(x)$, such that
%
%
\begin{equation}\label{eq:scb}
\lim_{n \to\infty} \mathsf{P}\bigl(\mu_n^{\mathrm{lower}}(x) \le\mu(x)
\le\mu_n^{\mathrm{upper}}(x) \mbox{ for all } x \in T\bigr) = 1-\alpha.
\end{equation}
The SCB can be used for model validation: one can test whether
$\mu(\cdot)$ is of certain parametric functional form by checking
whether the fitted parametric form lies in the SCB. Following the
work of \citet{BR73}, \citet{john82} derived the
asymptotic distribution of ${\sup_{0 \le x \le1}} |\mu_n(x) -
{\mathsf{E}}[\mu_n(x)]|$, assuming that $(X_i, Y_i)$ are independent random
samples from a bivariate population. Johnston's derivation is no
longer valid if dependence is present. For other work on
regression confidence bands under independence see \citet{KSY85},
\citet{HT88}, \citet{HM91}, \citet{SL94}, \citet{xia98},
\citet{CFN01} and \citet{dumb03}, among others. Recently
\citet{ZW08} proposed a method for constructing SCB for
stochastic regression models which have asymptotically correct
coverage probabilities. However, their confidence band is over an
increasingly dense grid of points instead of over an interval [see
also \citet{buhl98} and \citet{KSY85}]. Here we shall
also solve the latter problem and establish a uniform asymptotic
theory for the regression estimate $\mu_n(x)$, so that one can
construct a genuine SCB for regression functions. A similar result
will be derived for $\sigma(\cdot)$ as well.

The rest of the paper is organized as follows. Main results are
presented in Section \ref{sec:main}. Proofs are given in Sections
\ref{sec:proof} and \ref{sec:pfth45}. Our results are applied in
Section~\ref{sec:tb} to the U.S. Treasury yield rates data.

\section{Main results}
\label{sec:main} Before stating our theorems, we first introduce
dependence measures. Assume $X_k \in\mathcal{L}^p$, $p > 0$. Here
for a random variable $W$, we write $W \in\mathcal{L}^p$ ($p > 0$),
if $\| W \|_p := (\mathsf{E}|W|^p)^{1/p} < \infty$. Let
$\{\varepsilon'_j\}_{j \in{\mathsf{Z}}}$ be an i.i.d. copy of
$\{\varepsilon_j\}_{j \in{\mathsf{Z}}}$; let $\xi_n = (\ldots,
\varepsilon_{n-1}, \varepsilon_{n})$ and
\[
X'_n = G(\xi_n')\qquad\mbox{where }
\xi_n' = (\xi_{-1}, \varepsilon'_{0},
\varepsilon_1, \ldots,\varepsilon_n).
\]
Here $X_n'$ is a coupled process of $X_n$ with $\varepsilon_{0}$
in the latter replaced by an i.i.d. copy~$\varepsilon_0'$. Following
\citet{wu05}, define the physical dependence measure
\[
\theta_{n,p}=\|X_{n}-X'_{n}\|_{p}.
\]
Let $\theta_{n, p} = 0$ if $n < 0$. A similar quantity can be
defined if we couple the whole past: let $\xi^{\star}_{k, n} =
(\ldots, \varepsilon'_{k-n-2}, \varepsilon'_{k-n-1}, \xi_{k-n,k}),
k \ge n, $ where $\xi_{i, j} = (\varepsilon_i,
\varepsilon_{i+1}, \ldots, \varepsilon_j)$, and define
%
%
\begin{equation}\label{eq:J9-01}
\Psi_{n, p}=\|G(\xi_n)-G(\xi^{\star}_{n,n})\|_p.
\end{equation}
Our conditions on dependence will be expressed in terms of
$\theta_{n,p}$ and $\Psi_{n, p}$.

\subsection{Kernel density estimates}

We first consider a special case of (\ref{a1}) in which $X_n$ has
the form
%
%
\begin{equation}\label{eq:add}
X_{n} = a_{0}\varepsilon_{n} +
g(\ldots,\varepsilon_{n-2},\varepsilon_{n-1})
= a_{0}\varepsilon_{n} + g(\xi_{n-1}),
\end{equation}
where $g$ is a measurable function and $a_{0} \neq0$. Then the
coupled process $X'_n = a_0 \varepsilon_n + g(\xi_{-1},
\varepsilon'_0, \varepsilon_1, \ldots,\varepsilon_{n-1})$. We need
the following conditions:

(C1). There exists $0 < \delta_{2} \le\delta_{1} < 1$ such that
$n^{-\delta_{1}} = O(b_n)$ and $b_n = O( n^{-\delta_{2}})$.

(C2). Suppose that $X_1 \in\mathcal{L}^p$ for some $p > 0$. Let $p'
= \min(p, 2)$ and $\Theta_n = \sum_{i=0}^n \theta_{i,p'}^{p'/2}$.
Assume $\Psi_{n,p'}=O(n^{-\gamma})$ for some $\gamma> \delta_1
/(1-\delta_{1})$ and
%
%
\begin{equation}\label{eq:C275}
\mathcal{Z}_{n} b n^{-1} = o(\log n)\qquad\mbox{where }
\mathcal{Z}_{n}=
\sum_{k=-n}^{\infty}(\Theta_{n+k}-\Theta_{k})^{2}.
\end{equation}

(C3). The density function $f_{\varepsilon}$ of $\varepsilon_{1}$
is positive and
\[
\sup_{x\in{\mathsf{R}}}[f_{\varepsilon}(x)
+|f'_{\varepsilon}(x)|+|f''_{\varepsilon}(x)|]<\infty.
\]

(C4). The support of $K$ is $[-A, A]$, where $K$ is differentiable over
$(-A, A)$, the right (resp., left) derivative $K'(-A)$ [resp., $K'(A)$]
exists, and\break ${\sup_{|x| \le A}} |K'(x)| < \infty$. The Lebesgue measure
of the set $\{x\in[-A, A]\dvtx K(x)=0\}$ is zero. Let $\lambda_K = \int
K^{2}(y) \,d y$, $K_1 = [K^{2}(-A)+K^{2}(A)] / (2 \lambda_K)$ and $K_2 =
\int_{-A}^{A} (K'(t))^2 \,d t / (2 \lambda_K)$.
\begin{theorem}
\label{th:21} Let $l,u\in{\mathsf{R}}$ be fixed and $X_{n}$ be of form
(\ref{eq:add}). Assume \textup{(C1)--(C4)}. Then we have for every $z
\in
{\mathsf{R}}$,
%
%
\begin{equation}\label{th1}
\mathsf{P} \bigl((2\log\bar{b}{}^{-1})^{1/2} (\Delta
_{n}-d_{n} )
\leq z \bigr)\to e^{-2e^{-z}},
\end{equation}
where $\bar{b} = b / (u-l)$,
\[
d_{n}=(2\log\bar{b}^{-1})^{1/2}
+ \frac{1}{(2\log\bar{b}^{-1})^{1/2}}
\biggl\{ \log\frac{K_{1}}{\pi^{1/2}}
+\frac{1}{2}\log\log\bar{b}^{-1} \biggr\},
\]
if $K_{1}>0$, and otherwise
\[
d_{n}=(2\log\bar{b}^{-1})^{1/2}+\frac{1}{(2\log\bar{b}^{-1})^{1/2}}
\log\frac{K_2^{1/2}}{2^{1/2}\pi}.
\]
\end{theorem}

We now discuss conditions (C1)--(C4). The bandwidth condition (C1)
is fairly mild. In (C2), the quantity $\Theta_{n}$ measures the
cumulative dependence of $X_0, \ldots, X_n$ on $\varepsilon_0$,
and, with (C1), it gives sufficient dependence and bandwidth
conditions for the asymptotic Gumbel convergence (\ref{th1}). For
short-range dependent linear process $X_n = \sum_{j=0}^{\infty}
a_{j} \varepsilon_{n-j}$ with $\mathsf{E}\varepsilon_{1}=0$ and
$\mathsf{E}
\varepsilon^2_1 = 1$, (C2) is satisfied if $\sum_{j=0}^{\infty}
|a_j| < \infty$ and $\sum_{j=n}^\infty a^{2}_{j}= O(n^{-\gamma})$
for some $\gamma>2\delta_{1}/(1-\delta_{1})$. The latter condition
can be weaker than $\sum_{j=0}^{\infty}|a_{j}|<\infty$ if
$\delta_1<1/3$. Interestingly, (C2) also holds for some long-range
dependent processes; see Theorem \ref{th:lrd}. With (C3), it is
easily seen that $X_i$ does have a density. If (C3) is violated,
then $X_i$ may not have a density. For example, if $\varepsilon_i$
are i.i.d. Bernoulli with $\mathsf{P}( \varepsilon_i = 0) = \mathsf{P}(
\varepsilon_i = 1) = 1/2$, then $X_0 = \sum_{i=0}^\infty\rho^i
\varepsilon_{-i}$, where $\rho= (\sqrt5 - 1)/2$, does not have a
density [\citet{erdos39}]. The kernel condition (C4) is quite mild
and it is satisfied by many popular kernels. For example, it holds
for the Epanechnikov kernel $K(u) = 0.75 (1-u^2) \mathbf{1}_{|u| \le
1}$.

In Theorem \ref{th:ifs} below, we do not assume the special form
(\ref{eq:add}). We need regularity conditions on conditional
density functions. For jointly distributed random vectors $\xi$
and $\eta$, let $F_{\eta|\xi}(\cdot)$ be the conditional
distribution function of $\eta$ given $\xi$; let
$f_{\eta|\xi}(x)=\partial F_{\eta|\xi}(x)/\partial x$ be the
conditional density. For function $g$ with $\mathsf{E}|g(\eta)| <
\infty$, let $\mathsf{E}(g(\eta)|\xi)=\int g(x)\,d F_{\eta|\xi}(x)$
be the
conditional expectation of $g(\eta)$ given $\xi$.

Conditions (C2) and (C3) are replaced, respectively, by:

(C2)$'$. Suppose that $X_1 \in\mathcal{L}^p$ and $\theta_{n,p} =
O(\rho^{n})$ for some $p>0$ and $0<\rho<1$.

(C3)$'$. The density function $f$ is positive and there exists a
constant $B < \infty$ such that
\[
\sup_{x} [|f_{X_{n} |\xi_{n-1}} (x)|
+|f'_{X_{n} |\xi_{n-1}} (x)|
+|f''_{X_{n} |\xi_{n-1}} (x)|]\leq B \qquad\mbox{almost surely.}
\]
\begin{theorem}\label{th:22}
\label{th:ifs} Under \textup{(C1)}, \textup{(C2)$'$}, \textup{(C3)$'$}
and \textup{(C4)},
we have (\ref{th1}).
\end{theorem}

Many nonlinear time series models (e.g., ARCH models, bilinear
models, exponential AR models) satisfy (C2)$'$; see \citet{SW07}.
If $(X_i)$ is a Markov chain of the form $X_i = R(X_{i-1},
\varepsilon_i)$, where $R(\cdot, \cdot)$ is a bivariate measurable
function, then $f_{X_{i} |\xi_{i-1}} (\cdot)$ is the conditional
density of $X_i$ given $X_{i-1}$. Consider the ARCH model $X_i =
\varepsilon_i (a^2 + b^2 X_{i-1}^2)^{1/2}$, where $a > 0, b > 0$
are real parameters and $\varepsilon_i$ has density function
$f_\varepsilon$, then $f_{X_{i} |X_{i-1}} (x) = f_\varepsilon( x
/ H_i) / H_i$, where $H_i = (a^2 + b^2 X_{i-1}^2)^{1/2}$. So (C3)$'$
holds if $ \sup_{x} [f_\varepsilon(x) + |f'_\varepsilon(x)| +
|f''_\varepsilon(x)|] < \infty$ [cf. (C3)]. For more general
ARCH-type processes see \citet{DMR07}.

For short-range dependent processes for which
%
%
\begin{equation}\label{eq:C2751}
\Theta_\infty=\sum_{i=0}^\infty\theta_{i,p'}^{p'/2} < \infty,
\end{equation}
we have $\mathcal{Z}_{n} = O(n)$ and (\ref{eq:C275}) of condition
(C2) trivially holds. For long-range dependent processes,
(\ref{eq:C2751}) can be violated. A popular model for long-range
dependence is the fractionally integrated auto-regressive moving
average process [\citet{GJ80}, \citet{hosking81}]. Here
we consider the more general form of linear processes with slowly
decaying coefficients:
%
%
\begin{equation}\label{eq:lrdJ}
X_{n} = \sum_{j=0}^\infty a_j \varepsilon_{n-j}
\qquad\mbox{where } a_j = j^{-\beta} \ell(j), 1/2 < \beta< 1.
\end{equation}
Here $a_0 = 1$, $\ell(\cdot)$ is a slowly varying function and
$\varepsilon_i$ are i.i.d. with $\mathsf{E}\varepsilon_i = 0$ and
$\mathsf{E}
\varepsilon^2_i = 1$.
\begin{theorem}
\label{th:lrd} Assume (\ref{eq:lrdJ}). Let $l,u\in{\mathsf{R}}$ be fixed.
\textup{(i)} Assume \textup{(C1)}, \textup{(C3)}, \textup{(C4)},
$\delta_1 / (1-\delta_1) < \beta-1/2$ and
%
%
\begin{equation}\label{eq:Jsb}
b_n^{1/2} n^{1-\beta} \ell(n) = o(\log^{-1/2} n).
\end{equation}
Then (\ref{th1}) holds. \textup{(ii)} Assume \textup{(C1)}, \textup
{(C3)}, \textup{(C4)},
${\sup_x} |f'''_{\varepsilon}(x)| < \infty$ and
%
%
\begin{equation}\label{eq:Jlb}
\log^{1/2} n = o(b_n^{1/2} n^{1-\beta} \ell(n)).
\end{equation}
Let $c_\beta= \int_0^\infty(x+x^2)^{-\beta} \,d x / [(3-2\beta)
(1-\beta)]$. Then
%
%
\begin{equation}\label{eq:aug241}
{{\Delta_n} \over{b_n^{1/2} n^{1-\beta} \ell(n)}}
\Rightarrow|N(0,1)| { {\sqrt{c_\beta} \over\sqrt{\lambda_K}}}
\max_{l \le x \le u} { {|f'(x)| \over\sqrt{f(x)}}}.
\end{equation}
\end{theorem}

Theorem \ref{th:lrd} reveals the interesting dichotomy phenomenon
for the maximum deviation $\Delta_n$: if the bandwidth $b_n$ is
small such that (\ref{eq:Jsb}) holds, then the asymptotic
distribution is the same as the one under short-range dependence.
However, if $b_n$ is large, then both the normalizing constant and
the asymptotic distribution change. Let $b_n = n^{-\delta}
\ell_1(n)$, where $\ell_1$ is another slowly varying function.
Simple algebra shows that, if $\max((1+\delta) / (1-\delta), 2 -
\delta) < 2\beta$, then the bandwidth condition in Theorem
\ref{th:lrd}(i) holds. The latter inequality requires $\beta>
\sqrt3/2 = 0.866025,\ldots.$ If $\beta< 1 - \delta/ 2$, then
(\ref{eq:Jlb}) holds. Theorem \ref{th:lrd}(ii) is similar to
Theorem 3.1 in \citet{HH96}, with our result having a wider
range of $\beta$.

\subsection{Estimation of $\mu(\cdot)$ and $\sigma^2(\cdot)$}
Let $\widetilde{\xi}_i = (\ldots, \eta_{i-1}, \eta_i, \xi_i)$. For
a function $h$ with $\mathsf{E}h^2(\eta_i) < \infty$, write
\[
M^{r}_{n}(x)=\frac{1}{nb}
\sum_{k=1}^{n}K \biggl(\frac{X_{k}-x}{b} \biggr)Z_k\qquad
\mbox{where } Z_k = h(\eta_k) - \mathsf{E}h(\eta_k).
\]
\begin{proposition}
\label{th4} Let $l,u\in{\mathsf{R}}$ be fixed. Assume $\sigma
^{2}=\mathsf{E}
Z^{2}_{1}$ and $\mathsf{E}|Z_{1}|^{p} < \infty$, $p>2/(1-\delta_{1})$.
\textup{(i)} Assume (\ref{eq:add}), \textup{(C1)}, \textup{(C3)--(C4)} and
$\Psi_{n,q}=O(n^{-\gamma})$ for some $q>0$ and $\gamma>
\delta_{1}/(1-\delta_{1})$. Then for all $z \in{\mathsf{R}}$,
%
%
\begin{equation}\label{ps:1}
\mathsf{P} \Biggl( \sqrt{\frac{n b}{\lambda_K}}
\sup_{l\leq x\leq u} { {|M^{r}_{n}(x)|} \over f^{1/2}(x)\sigma}
- d_{n} \leq{z\over{(2\log\bar{b}^{-1})^{1/2}}} \Biggr)
\rightarrow e^{-2e^{-z}}
\end{equation}
as $n \to\infty$. \textup{(ii)} Assume\vspace*{1pt} (\ref{a1}), \textup
{(C1)}, \textup{(C2)$'$},
\textup{(C3)$'$} and \textup{(C4)} hold with $\xi_{n-1}$ in \textup
{(C2)$'$} replaced by
$\widetilde{\xi}_{n-1}$. Then (\ref{ps:1}) holds.
\end{proposition}

Proposition \ref{th4}(i) allows for long-range\vspace*{1pt} dependent
processes. For (\ref{eq:lrdJ}), by Karamata's theorem, $\Psi_{n,2} =
O(n^{1/2-\beta} \ell(n))$. So we have $\Psi_{n,2} = O(n^{-\gamma})$
with $\gamma> \delta_{1}/(1-\delta_{1})$ if $\delta_1 <
(2\beta-1)/(2\beta+1)$.

For $S \subset{\mathsf{R}}$, denote by $\mathcal{C}^p(S) =
\{g(\cdot)\dvtx
\sup_{x\in S}|g^{(k)}(x)|<\infty, k=0,\ldots,p\}$ the set of
functions having bounded derivatives on $S$ up to order $p\geq1$.
Let $S^{\epsilon} = \bigcup_{y\in S} \{x\dvtx|x-y| \le\epsilon\}$ be
the $\epsilon$-neighborhood of $S$, $\epsilon>0$.
\begin{theorem}
\label{th:scbreg} Let $l, u \in{\mathsf{R}}$ be fixed and $K$ be
symmetric. Assume that the conditions in Proposition \ref{th4}
hold with $Z_n = \eta_{n}$, $f_{\varepsilon}(\cdot), \mu(\cdot)
\in\mathcal{C}^{4} (T^{\epsilon})$ for some $\epsilon>0$, where
$T = [l, u]$, and that $b$ satisfies
\[
0<\delta_{1}<1/3,\qquad nb^{9}\log n=o(1)\quad
\mbox{and}\quad\mathcal{Z}_{n}b^{3}=o(n\log n).
\]
Let $\psi_K = \int u^2 K(u) \,d u / 2$ and $\rho_{\mu}(x) = \mu''(x)
+ 2 \mu'(x) f'(x) / f(x)$. Then
%
%
\begin{eqnarray}\label{eq:aug261}\quad
&&\mathsf{P} \Biggl( \sqrt{\frac{n b}{\lambda_K}}
\sup_{l\leq x\leq u}
{ {\sqrt{f_{n}(x)}|\mu_n(x)-\mu(x)-b^2 \psi_K \rho_{\mu}(x)|}
\over{\sigma(x)}} \nonumber\\[-8pt]\\[-8pt]
&&\hspace*{137.4pt}{}
- d_n \leq{z\over{(2\log\bar{b}^{-1})^{1/2}}} \Biggr)
\rightarrow e^{-2e^{-z}}.\nonumber
\end{eqnarray}
\end{theorem}

Note that $\sigma^2(x) = \mathsf{E}[ (Y_k - \mu(X_k))^2 | X_k = x
]$. It
is natural to use the Nadaraya--Watson method to estimate $\sigma^2(x)$
based on the residuals $\hat e_k = Y_{k}-\mu_n(X_{k})$:
\[
\sigma^{2}_{n}(x)=\frac{1}{n h f_{n1}(x)}
\sum_{k=1}^{n}K \biggl(\frac{X_{k}-x}{h} \biggr)[Y_{k}-\mu_n(X_{k})]^{2},
\]
where the bandwidths $h = h_n \to0$ and $n h_n \to\infty$, and
\[
f_{n1}(x)=\frac{1}{nh}
\sum_{k=1}^{n} K \biggl(\frac{X_{k}-x}{h} \biggr).
\]
\begin{theorem}
\label{th:scbreg-2} Let $l,u\in{\mathsf{R}}$ be fixed and $K$ be
symmetric. Assume $\nu_{\eta} =\mathsf{E}\eta^{4}_{1}-1 < \infty$.
Further assume that the conditions in Proposition \ref{th4} hold
with $Z_{n} = \eta^{2}_{n}-1$, $f(\cdot), \sigma(\cdot) \in
\mathcal{C}^{4} (T^{\epsilon})$ for some $\epsilon>0$, where $T =
[l, u]$, and that $h\asymp b$ satisfies
\[
0<\delta_{1}<1/4,\qquad nb^{9}\log n=o(1)
\]
and
\[
\mathcal{Z}_{n}b^{3}=o(n\log n).
\]
Let $\rho_{\sigma}(x) = 2{\sigma'}^2(x)+2\sigma(x)\sigma''(x) + 4
\sigma(x)\sigma'(x) f'(x) / f(x)$. Then
%
%
\begin{eqnarray}\label{eq:Aug262}\quad\qquad
&& \mathsf{P} \Biggl( \sqrt{\frac{n h}{\lambda_K \nu
_{\eta}}}
\sup_{l\leq x\leq u}
{ {\sqrt{f_{n1}(x)}|\sigma^{2}_n(x)-\sigma^{2}(x)
-h^2 \psi_K \rho_{\sigma}(x)|}
\over{\sigma^{2}(x)}} \nonumber\\[-8pt]\\[-8pt]
&&\hspace*{137.4pt}\hspace*{17.8pt}{}
- d_n \leq{z\over{(2\log\bar{h}^{-1})^{1/2}}} \Biggr)
\rightarrow e^{-2e^{-z}},\nonumber
\end{eqnarray}
where $d_{n}$ is defined as in Theorem \ref{th:21} by replacing
$\bar{b}$ with $\bar{h}=h/(u-l)$.
\end{theorem}

We now compare the SCBs constructed based on Theorem 1 in
\citet{ZW08} and Theorem \ref{th:scbreg}. Assume $l=0$ and $u=1$. The
former is over the grid point $T_n = \{ 2 b_n j, j=0, 1, \ldots,
J_n\}$ with $J_n = \lceil1/(2 b_n) \rceil$, while the latter is a
genuine SCB in the sense that it is over the whole interval $T =
[0, 1]$. Let $\hat\rho_{\mu}(\cdot)$ [resp., $\hat\sigma(\cdot)$] be a
consistent estimate of $\rho_{\mu}(\cdot)$ [resp., $\sigma(\cdot)$] and
$z_\alpha= -\log\log(1-\alpha)^{-1/2}$, $0 < \alpha< 1$. By
Theorem \ref{th:scbreg}, we can construct the $1-\alpha$ SCB for
$\mu(x)$ over $x \in[0, 1]$ as
%
%
\begin{eqnarray}\label{eq:aug6-1}
&&\mu_n(x) - b^{2} \psi_K \hat\rho_{\mu}(x)
\pm l_1 \hat\sigma(x) \sqrt{\frac{\lambda_K}{n b
f_n(x)}}\nonumber\\[-8pt]\\[-8pt]
&&\eqntext{\mbox{where }
\displaystyle l_1 = { {z_\alpha} \over{(2\log b^{-1})^{1/2}}} + d_n.}
\end{eqnarray}
Similarly, using Theorem 1 in \citet{ZW08}, the $1-\alpha$
confidence band for $\mu(x)$ over $x \in T_n$ is also of form
(\ref{eq:aug6-1}) with $l_1$ replaced by
\[
l_2 = { {z_\alpha} \over{(2\log J_n)^{1/2}}} + (2 \log J_n)^{1/2}
- { {{1/2} \log\log J_n + \log(2 \sqrt\pi)}
\over{(2 \log J_n)^{1/2}}}.
\]
Elementary calculations show that, interestingly, $l_1$ and $l_2$
are quite close: $l_1 - l_2 = (\log\log b^{-1}) / (2 \log
b^{-1})^{1/2} ( 1 + o(1))$ if $K_1 > 0$.

\section{Application to the treasury bill data}
\label{sec:tb} There is a huge literature on models for short-term
interest rates. Let $R_t$ be the interest rate at time $t$. Assume
that $R_t$ follows the diffusion model
%
%
\begin{equation}\label{eq:diffu}
d R_t = \mu(R_t) \,d t + \sigma(R_t) \,d \mathbb{B}(t),
\end{equation}
where $\mathbb{B}$ is the standard Brownian motion, $\mu(\cdot)$ is
the instantaneous return or drift function and $\sigma(\cdot)$ is
the volatility function. \citet{BS73} considered the
model with $\mu(x) = \alpha x$ and $\sigma(x) = \sigma x$. Vasicek
(\citeyear{V77}) assumed that $\mu(x) = \alpha_0 + \alpha_1 x$ and
$\sigma(x) \equiv\sigma$, where $\alpha_0, \alpha_1$ and $\sigma$
are unknown constants. \citet{CIR85} and
\citet{court82} assumed that $\sigma(x) = \sigma x^{1/2}$ and
$\sigma(x) = \sigma x$, respectively. Both models are generalized
by \citet{chanetal92} to the form $\sigma(x) = \sigma x^\gamma$,
with $\sigma$ and $\gamma$ being unknown parameters.
\citet{stant97}, \citet{fan98}, Chapman and Pearson (\citeyear
{ChP00}) and Fan and
Zhang (\citeyear{FZ03}) considered the nonparametric estimation of
$\mu(\cdot)$ and $\sigma(\cdot)$ in (\ref{eq:diffu}); see also
A\"{i}t-Sahalia (\citeyear{Ait1}, \citeyear{Ait2}). \citet
{stant97} constructed
\textit{point-wise} confidence intervals which serve as a tool for suggesting
which parametric models to use. \citet{zhao08} gave an excellent
review of parametric and nonparametric approaches of
(\ref{eq:diffu}). See also the latter paper for further
references.

%
\begin{figure}[b]

\includegraphics{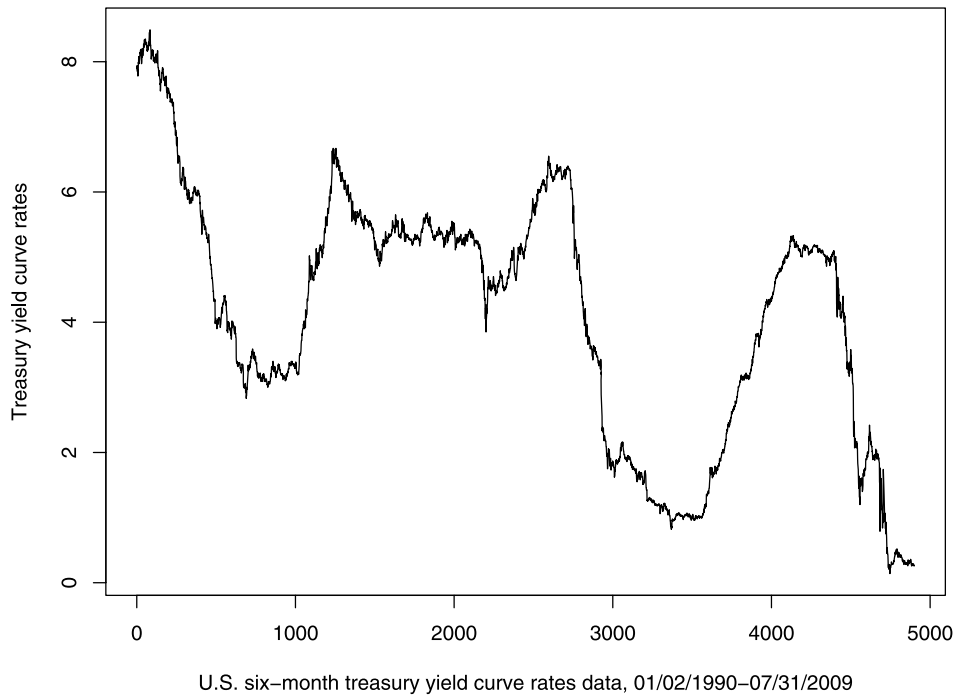}

\caption{U.S. six-month treasury yield curve rates data from
January 2nd, 1990 to July 31st, 2009. Source: U.S. Treasury
department's website \protect\url{http://www.ustreas.gov/}.}
\label{fid:6m}
\end{figure}

%
\begin{figure}

\includegraphics{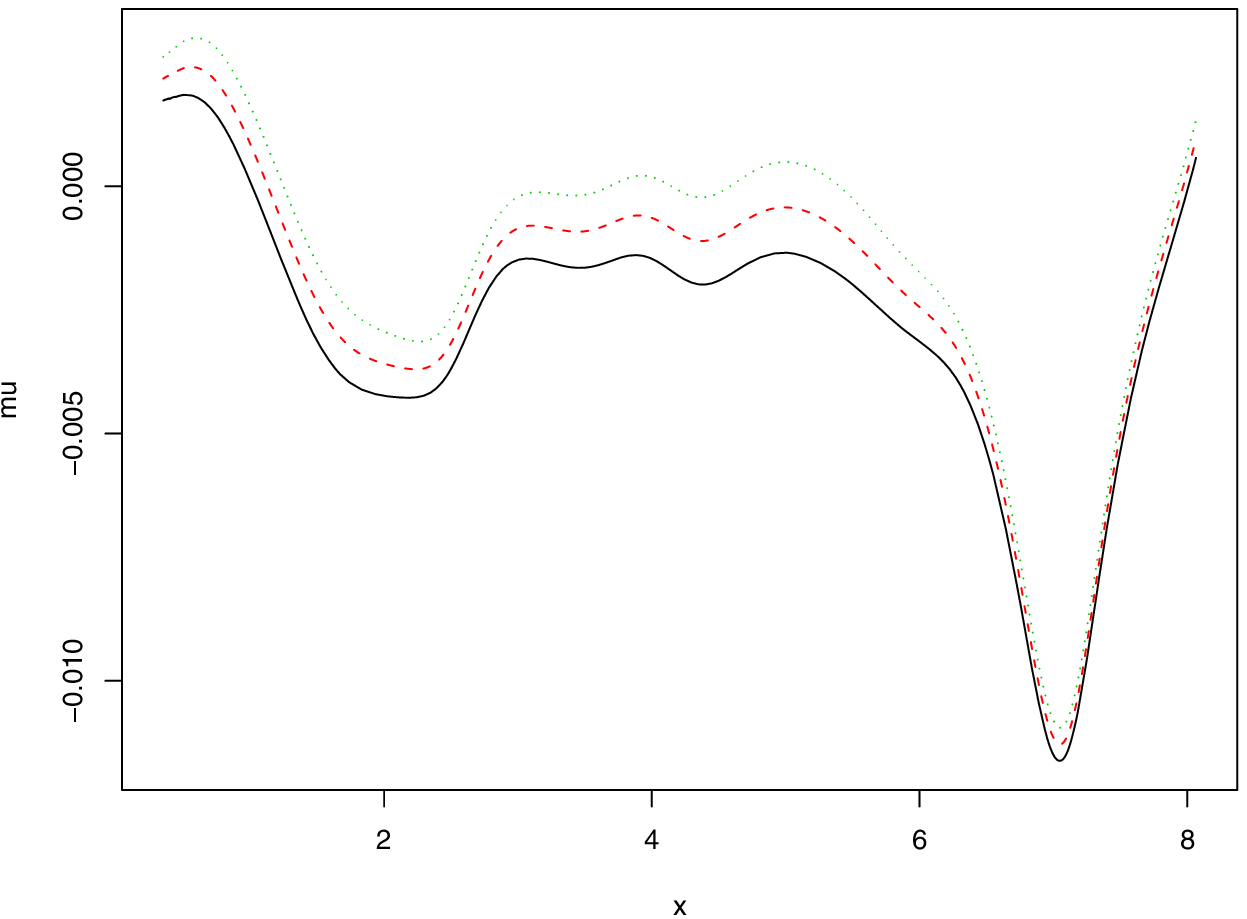}

\caption{$95\%$ SCB of the regression function $\mu(\cdot)$ over
the interval $[l, u] = [0.35, 8.06]$. The dashed curve in the
middle is $\mu_n(x)- b^{2} \psi_K \hat\rho(x)$, the
bias-corrected estimate of $\mu$. }\label{fid:scbmu}
\end{figure}

Here we shall consider the U.S. six-month treasury yield rates
data from January 2nd, 1990 to July 31st, 2009. The data can be
downloaded from the U.S. Treasury department's website
\url{http://www.ustreas.gov/}. It has 4900 daily rates and a plot is
given in Figure \ref{fid:6m}. Let $X_i = R_{t_i}$ be the rate at
day $i=1, \ldots, 4900$. For the daily data, since one year has
250 transaction days, $t_i - t_{i-1} = 1/250$. Let $\Delta=
1/250$. As a discretized version of (\ref{eq:diffu}), we consider
the model
%
%
\begin{equation}\label{eq:dsdiffu}
Y_i = \mu(X_i) \Delta+ \sigma(X_i) \Delta^{1/2} \eta_i,
\end{equation}
where $Y_i = R_{t_{i+1}} - R_{t_i} = X_{i+1} - X_i$ and $\eta_i =
(\mathbb{B}(t_{i+1}) - \mathbb{B}(t_i)) / \Delta^{1/2}$ are i.i.d. standard
normal. For convenience of applying Theorem \ref{th:scbreg}, in
the sequel we shall write $\mu(X_i) \Delta$ [resp., $\sigma(X_i)
\Delta^{1/2}$] in (\ref{eq:dsdiffu}) as $\mu(X_i)$ [resp.,
$\sigma(X_i)$]. So (\ref{eq:dsdiffu}) is rewritten as
%
%
\begin{equation}\label{eq:dsdiffu2}
Y_i = \mu(X_i) + \sigma(X_i) \eta_i.
\end{equation}

Figure \ref{fid:scbmu} shows the estimated $95\%$ simultaneous
confidence band for the regression function $\mu(\cdot)$ over the
interval $T = [l, u] = [0.35, 8.06]$, which includes $96 \%$ of
the daily rates $X_i$. To select the bandwidth, we use the \texttt{R}
program \texttt{bw.nrd} which gives $b = 0.37$. Then we use the \texttt
{R} program \texttt{locpoly} for local polynomial regression. The
Nadaraya--Watson estimate is a special case of the local polynomial
regression with degree $0$. The function $\rho(x)$ in the bias
term $b^2 \psi_K \rho(x)$ in Theorem \ref{th:scbreg} involves the
first and second order derivatives $\mu'$, $f'$ and $\mu''$. The
program \texttt{locpoly} can also be used to estimate derivatives
$\mu'$ and $\mu''$, where we use the bigger bandwidth $2 b =
0.74$. For $f$, we use the \texttt{R} program \texttt{density}, and
estimate $f'$ by differentiating the estimated density. Then we
can have the bias-corrected estimate $\tilde\mu_n(x) = \mu_n(x) -
b^{2} \psi_K \hat\rho(x)$ for $\mu$, which is plotted in the the
middle curve in Figure \ref{fid:scbmu}. To estimate
$\sigma(\cdot)$, as in \citet{stant97}, we shall make use of the
estimated residuals $\hat e_i = Y_i - \tilde\mu_n(X_i)$, and
perform the Nadaraya--Watson regression of $\hat e_i^2$ versus
$X_i$ with the bandwidth $b$. In our data analysis the boundary
problem of the Nadaraya--Watson regression raised in Chapman and
Pearson (\citeyear{ChP00}) is not severe since we focus on the interval
$T =
[0.35, 8.06]$, while the whole range is $[\min X_i, \max X_i] =
[0.14, 8.49]$.

The Gumbel convergence in Theorem \ref{th:scbreg} can be quite
slow, so the SCB in (\ref{eq:aug6-1}) may not have a good
finite-sample performance. To circumvent this problem, we shall
adopt a simulation based method. Let
\[
\Pi_n = \sup_{x\in T} {
{|\sum_{k=1}^{n}K(X^*_{k}/b-x/b) \eta^*_{k}|}
\over{n b f^{1/2}(x)} },
\]
where $X^*_k$ are i.i.d. with density $f$, $\eta^*_{k}$ are i.i.d. with
$\mathsf{E}\eta_{n}=0$, $\mathsf{E}\eta^2_{n}=1$ and $\mathsf
{E}|\eta_{1}|^{p} <
\infty$, and $(X^*_k)$ and $(\eta^*_{k})$ are independent. As in
Theorem \ref{th:scbreg}, let
\[
\Pi'_n = \sup_{x\in T}
{ {\sqrt{f(x)}|\mu_n(x)-\mu(x)-b^2 \psi_K \rho(x)|}
\over{\sigma(x)}}.
\]
By Theorem \ref{th:scbreg} and Proposition \ref{th4}, with proper
centering and scaling, $\Pi_n$ and $\Pi_n'$ have the same
asymptotic Gumbel distribution. So the cutoff value, the
$(1-\alpha)$th quantile of $\Pi'_n$, can be estimated by the
sample $(1-\alpha)$th quantile of many simulated $\Pi_n$'s. For the
U.S. Treasury bill data, we simulated $10\mbox{,}000$ $\Pi_n$'s and obtained
the $95\%$ sample quantile $0.39$. Then the SCB is constructed as
$\tilde\mu_n(x) \pm0.39 \hat\sigma(x) / f^{1/2}_n(x)$; see the
upper and lower curves in Figure \ref{fid:scbmu}.

\begin{figure}

\includegraphics{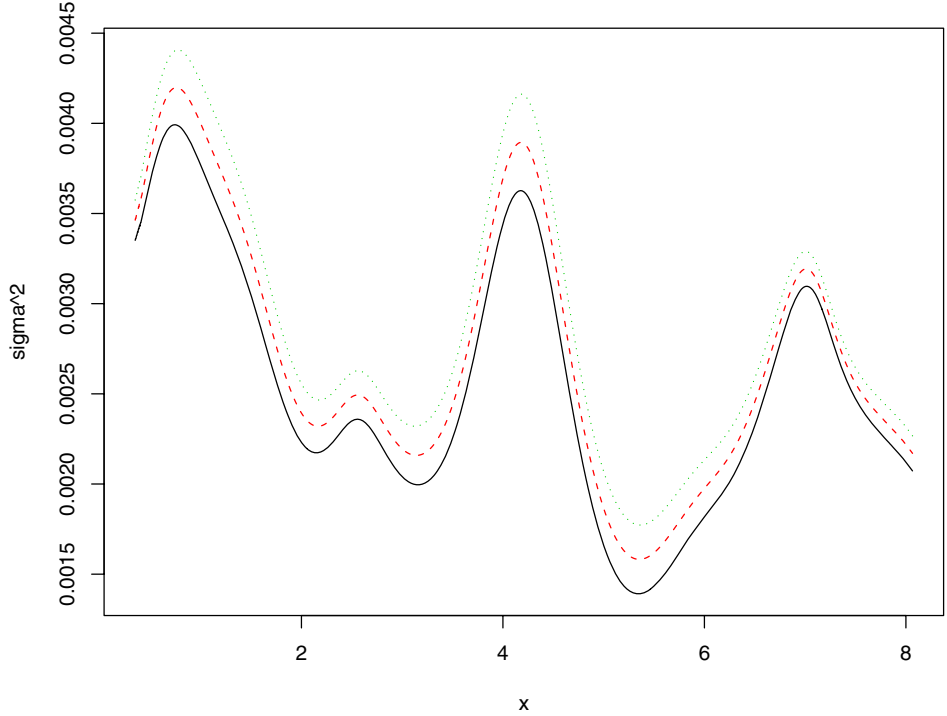}

\caption{$95\%$ SCB of the volatility function $\sigma^2(\cdot)$
over the interval $[l, u] = [0.35, 8.06]$. The dashed curve in the
middle is $\sigma^2_n(x) - b^{2} \psi_K \hat\rho_\sigma(x)$, the
bias-corrected estimate of $\sigma^2$. }\label{fid:scbs2}
\end{figure}

We now apply Theorem \ref{th:scbreg-2} to construct SCB for
$\sigma^2(\cdot)$. We choose $h = b$, which has a reasonably
satisfactory performance in our data analysis. By Theorem
\ref{th:scbreg-2},
\[
\Pi''_n = {1\over\sqrt\nu_{\eta}} \sup_{x\in T}
{ {\sqrt{f(x)}|\sigma^{2}_n(x)-\sigma^{2}(x)
-b^2 \psi_K \rho_{\sigma}(x)|}
\over{\sigma^{2}(x)}}
\]
has the same asymptotic distribution as $\Pi_n$ and $\Pi'_n$.
Based on the above simulation, we choose the cutoff value $0.39$.
As in the treatment of $\mu'$ and $\mu''$ in the bias term of
$\mu_n$, we use a similar estimate, noting that $\rho_{\sigma}(x)
= (\sigma^2(x))''+2(\sigma^2(x))' f'(x) / f(x)$ has the same form
as $\rho_\mu(x)$. The $95\%$ SCB of $\sigma^2(\cdot)$ is presented
in Figure \ref{fid:scbs2}.

%

Based on the $95\%$ SCB of $\mu(\cdot)$, we conclude that the
linear drift function hypothesis $H_0\dvtx \mu(x) = \alpha_0 +
\alpha_1 x$ for some $\alpha_0$ and $\alpha_1$ is rejected at the
$5\%$ level. Other simple parametric forms do not seem to exist.
Similar claims can be made for $\sigma^2(\cdot)$, and none of the
parametric forms previously mentioned seems appropriate. This
suggests that the dynamics of the treasury yield rates might be
far more complicated than previously speculated.

\section{\texorpdfstring{Proofs of Theorems
\protect\ref{th:21}--\protect\ref{th:lrd}}{Proofs of Theorems 2.1--2.3}}
\label{sec:proof}
Throughout the proofs
$C$ denotes constants which do not depend on $n$ and $b_n$. The
values of $C$ may vary from place to place. Let $\lfloor\cdot
\rfloor$ and $\lceil\cdot\rceil$ be the floor and ceiling
functions, respectively. Without loss of generality, we assume
$l=0$, $u=1$ in (\ref{eq:J91}) and $A=1$ in condition (C4). Write
\[
{\sqrt{n b} \over{\sqrt{\lambda_K f(b t)}}}
[f_n(b t) - \mathsf{E}f_n(b t)] =
M_{n}(t) + N_{n}(t),
\]
where $M_n(t)$ has summands of martingale differences
\[
M_{n}(t)=\frac{1}{\sqrt{nb\lambda_K f(bt)}}\sum_{k=1}^{n}
\{K (X_{k}/b-t )-\mathsf{E}[ K ( X_{k}/ b-t ) |\xi_{k-1}] \},
\]
and, since $\mathsf{E}[ K (X_{k}/b - t) | \xi_{k-1} ] = b \int_{-1}^1
K(v) f_{X_{k}|\xi_{k-1}}(b v + b t) \,d v$, the remainder
\begin{eqnarray*}
N_{n}(t)&=&\frac{1}{\sqrt{nb\lambda_K f(bt)}}
\sum_{k=1}^{n} \{\mathsf{E}[ K ( X_{k}/ b-t )|\xi_{k-1}]
- \mathsf{E}K ( X_{k}/ b-t ) \} \\
&=& {\sqrt{b} \over{\sqrt{n \lambda_K f(b t)}}}
\int_{-1}^1 K(v) Q'_n(bv+bt) \,d v,
\end{eqnarray*}
where
\[
Q_n(x) = \sum_{k=1}^n [F_{X_{k}|\xi_{k-1}}(x)-F(x)].
\]
If $X_{n}$ admits the form (\ref{eq:add}), we assume $a_{0}=1$.
Let $Y_k = g(\ldots, \varepsilon_{k-1}, \varepsilon_k)$. Then
$f_{X_{k}|\xi_{k-1}}(bv+bt)=f_\varepsilon(b v + b t-Y_{k-1})$.
\begin{pf*}{Proofs of Theorems \ref{th:21} and \ref{th:22}}
We split $[1,n]$ into alternating big and small blocks $H_{1}$,
$I_{1}, \ldots, H_{\iota_n}$, $I_{\iota_n}$, $I_{\iota_n+1}$,
with length $|H_{i}| = \lfloor n^{\tau_{1}}\rfloor$, $|I_{i}| =
\lfloor n^{\tau}\rfloor$, $1\leq i\leq\iota_n$,
$|I_{\iota_n+1}|=n-\iota_n(\lfloor n^{\tau_{1}}\rfloor+\lfloor
n^{\tau}\rfloor)$ and $\iota_n = \lfloor n/(\lfloor
n^{\tau_{1}}\rfloor+\lfloor n^{\tau}\rfloor) \rfloor$, where
$\delta_1 / \gamma< \tau<\tau_{1} < 1-\delta_1$. Let $m=|I_{1}|$,
\begin{eqnarray*}
u_{j}(t)&=&\sum_{k\in H_{j}}
\{\mathsf{E}[ K ( X_{k}/b-t ) | \xi_{k-m,k} ]
-\mathsf{E}[ K ( X_{k}/b-t ) | \xi_{k-m,k-1} ]\},\\
v_{j}(t)&=&\sum_{k\in I_{j}} \{\mathsf{E}[K ( X_{k}/ b-t ) | \xi_{k-m,k}]
-\mathsf{E}[ K ( X_{k}/b-t ) | \xi_{k-m,k-1}] \},\\
\widetilde{M}_{n}(t)&=&\frac{1}{\sqrt{nb\lambda_K f(bt)}}
\sum_{j=1}^{\iota_n}u_{j}(t),\qquad
R_{n}(t)=\frac{1}{\sqrt{nb\lambda_K f(bt)}}
\sum_{j=1}^{\iota_n+1}v_{j}(t).
\end{eqnarray*}
Theorems \ref{th:21} and \ref{th:22} follow from Lemmas
\ref{le1}--\ref{le3} and Lemma \ref{lem:J9-1} below.
\end{pf*}
\begin{pf*}{Proof of Theorem \ref{th:lrd}}
Case (i) follows from Theorem \ref{th:21}. For (ii), since
$\sum_{i=1}^n Y_{i-1} / (c_\beta n^{3/2-\beta} \ell(n))
\Rightarrow N(0,1)$ [cf. \citet{HH96}],\vspace*{1pt} where
$Y_{i-1}=\sum_{k=1}^{\infty}a_{k}\varepsilon_{i-k}$, it follows
from (\ref{eq:Jlb}), Lemma \ref{le1}(ii) and Lemma
\ref{lem:Aug202}.
\end{pf*}
\begin{lemma}
\label{le1} Assume \textup{(C4)}. \textup{(i)} We have
%
%
\begin{equation}\label{eq:J12-2}
{\sup_{0\leq t\leq b^{-1}}}|N_{n}(t)|
= O_\mathsf{P}(b^{1/2} n^{-1/2} \widetilde{\Theta}_n),
\end{equation}
where $\widetilde{\Theta}_n=\mathcal{Z}^{1/2}_{n}$ if $(X_{n})$
satisfies (\ref{eq:add}) and \textup{(C3)}; $\widetilde{\Theta}_n =
O(n^{1/2})$ if $(X_{n})$ satisfies (\ref{a1}), \textup{(C2)$'$} and
\textup{(C3)$'$}. \textup{(ii)} For the process (\ref{eq:lrdJ}), we have
(\ref{eq:J12-2}) with $\widetilde\Theta_n = O(n^{3/2-\beta}
\ell(n))$, and
%
%
\begin{equation}\label{eq:J12-3}\quad
\sup_{0\leq t\leq b^{-1}}\Biggl|N_{n}(t) {\sqrt{n b \lambda_K f(b t)}}
- b f'(b t) \sum_{j=1}^n Y_{j-1}\Biggr|
= o(b n^{3/2-\beta} \ell(n)),
\end{equation}
where $Y_{j-1}=\sum_{k=1}^{\infty}a_{k}\varepsilon_{j-k}$.
\end{lemma}
\begin{lemma}
\label{le2} Under conditions of Theorems \ref{th:21} or
\ref{th:22}, we have
\[
\mathsf{P} \Bigl({\sup_{0\leq t\leq
b^{-1}}}|M_{n}(t)-\widetilde{M}_{n}(t)-R_{n}(t)| \geq(\log
b^{-1})^{-2} \Bigr)=o(1).
\]
\end{lemma}
\begin{lemma}
\label{le3} Under conditions of Theorems \ref{th:21} or
\ref{th:22}, we have
%
%
\begin{equation}\label{eq:AugR1}
\mathsf{P} \Bigl({\sup_{0\leq t\leq b^{-1}}}|R_{n}(t)|
\geq(\log b^{-1})^{-2} \Bigr)=o(1).
\end{equation}
\end{lemma}
\begin{lemma}\label{lem:Aug202}
Let $\sup_{x} f_{X_{n} |\xi_{n-1}}(x)$ be a.s. bounded. Assume
\textup{(C4)}. Then
\[
{\sup_{0 \le t \le b^{-1}} }|M_{n}(t)| = O_\mathsf{P}\bigl(\sqrt{\log
n} \bigr).
\]
Consequently, under conditions of Lemma \ref{le1}, $\mathsf{E}f_n(x) -
f(x) = f''(x) b^2 \psi_K + o(b^2)$ and
\[
{\sup_{0 \le x \le1} }|f_{n}(x)-f(x)|
= { {O_\mathsf{P}(\sqrt{\log n})} \over\sqrt{n b}}
+ { {O_\mathsf{P}(\tilde\Theta_n)} \over n} + O(b^2).
\]
\end{lemma}

Lemma \ref{lem:Aug202} gives an upper bound of $\sup_{0 \le t \le
b^{-1}} |M_{n}(t)|$. Under stronger conditions, one can have a far
deeper asymptotic distributional result. By Lemmas \ref{lem:J9-1},
\ref{le2} and \ref{eq:AugR1}, it is asymptotically distributed as
Gumbel.
\begin{lemma}
\label{lem:J9-1} Under conditions of Theorems \ref{th:21} or
\ref{th:22}, we have for all $z \in{\mathsf{R}}$ that
%
%
\begin{equation}\label{eq3}
\mathsf{P} \Bigl({\sup_{0\leq t\leq b^{-1}}}|\widetilde
{M}_{n}(t)|<x_z \Bigr)
\to e^{-2e^{-z}}\qquad\mbox{where }
x_z = d_{n} + {z\over{(2\log b^{-1})^{1/2}}}.\hspace*{-22pt}
\end{equation}
\end{lemma}

\subsection{\texorpdfstring{Proofs of Lemmas
\protect\ref{le1}--\protect\ref{lem:Aug202}}{Proofs of Lemmas 4.1--4.4}}

\mbox{}

\begin{pf*}{Proof of Lemma \ref{le1}}
We claim that, for any $a_{0}>0$,
%
%
\begin{equation}\label{eq1-1}
\mathsf{E} \Bigl[ {\sup_{|x|\leq a_{0}}}|Q'_n(x)|^2 \Bigr]
= O(\widetilde{\Theta}^{2}_n),
\end{equation}
which implies Lemma \ref{le1}(i) in view of
%
%
\begin{equation}\label{eq:aug201}
N_{n}(t)=\frac{\sqrt{b}}{\sqrt{n \lambda_K f(b t)}}
\int_{-1}^{1} K(x) Q'_{n}\bigl(b(x+t)\bigr) \,d x
\end{equation}
by noting that $\inf_{0 \le x \le1} f(x) > 0$, ${\int_{-1}^1}
|K(u)|\,du < \infty$. To prove (\ref{eq1-1}), we use Lemma 4 in
\citet{wu03}, which implies that
\[
{\sup_{|x|\leq a_{0}}}|Q'_n(x)|^2
\leq
2 a^{-1}_{0} \int_{-a_{0}}^{a_{0}}|Q'_n(x)|^2 \,d x
+ 2a_{0}\int_{-a_{0}}^{a_{0}}|Q''_n(x)|^2 \,d x.
\]
We first suppose that $(X_{n})$ satisfies (\ref{eq:add}) and (C3).
Let
\[
\mathcal{P}_{k} \cdot=\mathsf{E}( \cdot|\mathcal{F}_{k}) -\mathsf
{E}( \cdot
|\mathcal{F}_{k-1}), \qquad k \in{\mathsf{Z}},
\]
be the projection operators. By the orthogonality of
$\mathcal{P}_{k}$, we have
\begin{eqnarray*}
\|Q'_{n}(x)\|^{2}_{2}
&=&\sum_{k=-\infty}^{n}\|\mathcal{P}_{k}Q'_{n}(x)\|^{2}_{2}
\leq\sum_{k=-\infty}^{n} \Biggl(\sum_{i=1}^{n} \|\mathcal{P}_{k}
f_{X_{i}|\xi_{i-1}}(x)\|_{2} \Biggr)^{2}\\
&\leq&C\sum_{k=-\infty}^{n} \Biggl(\sum_{i=1-k}^{n-k}
\theta^{p'/2}_{i,p'} \Biggr)^{2} = C\mathcal{Z}_n,
\end{eqnarray*}
where $C$ does not depend on $x$. Similarly, we have $\sup_{x\in
{\mathsf{R}}}
\|Q''_{n}(x)\|^{2}_2 \leq C\mathcal{Z}_n$. This proves (\ref{eq1-1}).

To prove (\ref{eq1-1}) for $(X_{n})$ satisfying (\ref{a1}), (C2)$'$
and (C3)$'$, we note that
\begin{eqnarray*}
\sup_{x\in{\mathsf{R}}}\|\mathcal{P}_{k} F_{X_{i}|\xi_{i-1}}(x)\|^{2}_{2}
&\leq& \sup_{x\in{\mathsf{R}}}\mathsf{E}\bigl|I\{X_{i}\leq x\}-I\bigl\{
X_{i,\{k\}}\leq
x\bigr\}\bigr|
\\
&\leq& \sup_{x\in{\mathsf{R}}}
\mathsf{P} \bigl(|X_{i}-x|\leq\bigl|X_{i}-X_{i,\{k\}}\bigr| \bigr)\\
&\leq& C(\theta^{1/2}_{i-k,p}+\theta^{p/2}_{i-k,p}),
\end{eqnarray*}
where $X_{i,\{k\}}=G(\xi_{k-1},\varepsilon'_{k},\xi_{k+1,i})$ and
we used the inequality
\[
|I\{X\leq x\}-I\{Y\leq x\}|\leq I\{|X-x|\leq|X-Y|\}.
\]
Since ${\sup_{x}} |f'_{X_{n} |\xi_{n-1}}
(x)|\leq B$, we have
\[
\biggl|f_{X_{i}|\xi_{i-1}}(x)
-\frac{F_{X_{i}|\xi_{i-1}}(x)-F_{X_{i}|\xi_{i-1}}(x-\Delta)}
{\Delta} \biggr| \leq B\Delta,
\]
which by letting $\Delta= (\theta^{1/2}_{i-k,p} +
\theta^{p/2}_{i-k, p})^{1/2}$ yields that
\[
\sup_{x\in{\mathsf{R}}}\|\mathcal{P}_{k} f_{X_{i}|\xi_{i-1}}(x)\|^{2}_{2}
\leq C(\theta^{1/2}_{i-k,p}+\theta^{p/2}_{i-k,p})^{1/2}.
\]
This implies $\sup_{x\in{\mathsf{R}}}\|Q'_{n}(x)\|^{2}_{2}=O(n)$.
Similarly, we have $\sup_{x\in{\mathsf{R}}}\|Q''_{n}(x)\|
^{2}_{2}=O(n)$. We
finish the proof of Lemma \ref{le1}(i).

We now prove (\ref{eq:J12-3}). For $i \ge2$ write $Y_{i-1} = U +
a_i \varepsilon_0 + W$, where $U = \sum_{j=1}^{i-1} a_j
\varepsilon_{i-j}$ and $W = \sum_{j=i+1}^\infty a_j
\varepsilon_{i-j}$. Let $W' = \sum_{j=i+1}^\infty a_j
\varepsilon'_{i-j}$. Let $c_0 = \sup_x [|f_\varepsilon'(x)| +
|f_\varepsilon''(x)|]$. By Taylor's expansion, there exists $R \in
[0, 1]$ such that
\begin{eqnarray*}
\vartheta_i :\!&=& \sup_x \| f_\varepsilon(x-Y_{i-1}) -
f_\varepsilon(x-U - W)
+ a_i\varepsilon_0 f'_\varepsilon(x-U-a_i\varepsilon_0'-W')\|\\
&=& \sup_x \|{-}
a_i\varepsilon_0 f'_\varepsilon(x-U-R a_i\varepsilon_0-W)
+ a_i\varepsilon_0 f'_\varepsilon(x-U-a_i\varepsilon_0'-W')\|\\
&\le& \bigl\| a_i\varepsilon_0 c_0
\min(1, |a_i\varepsilon_0'|+ |a_i\varepsilon_0|+ |W| + |W'|)
\bigr\| = o(|a_i|).
\end{eqnarray*}
Here we use the fact that $\| \varepsilon_0 \min(1,
|a_i\varepsilon_0|) \| \to0$ since $a_i \to0$, and $a_i
\varepsilon_0$ and $|W| + |W'|$ are independent. Since
$\varepsilon'_l, \varepsilon_m$, $l, m \in{\mathsf{Z}}$, are i.i.d.,
we have
$f(x) = \mathsf{E}[ f_\varepsilon(x - U - a_i \varepsilon_0' - W') |
\xi_0]$. By the Lebesgue dominated convergence theorem, $f'(x) =
\mathsf{E}[ f'_\varepsilon(x - U - a_i \varepsilon_0' - W') | \xi
_0]$. By
Jensen's inequality,
\[
{\sup_x} \| \mathsf{E}[ f_\varepsilon(x-Y_{i-1}) -
f_\varepsilon(x-U - W) | \xi_0]
+ a_i\varepsilon_0 f'(x)\| \le\vartheta_i,
\]
which again by Jensen's inequality implies that $\sup_x \| \mathsf{E}[
f_\varepsilon(x-Y_{i-1}) - f_\varepsilon(x-U - W) | \xi_{-1}] \le
\vartheta_i$. Since $\mathsf{E}[ f_\varepsilon(x-U - W) | \xi_{-1}] =
\mathsf{E}[ f_\varepsilon(x-U - W) | \xi_0]$, we have
\[
{\sup_x} \| \mathcal{P}_0 [f_\varepsilon(x-Y_{i-1})
+ f'(x) Y_{i-1}] \| \le2 \vartheta_i = o(|a_i|).
\]
Define $\vartheta_i = 0$ if $i < 0$. Let $T_n(x) = Q_n(x) + f(x)
\sum_{i=1}^n Y_{i-1}$. If $k \le-n$, then
\[
\|\mathcal{P}_k T'_n(x)\|
\le\sum_{j=1}^n 2 \vartheta_{j-k}
= o(n |k|^{-\beta} \ell(|k|)).
\]
If $-n < k \le n$, by Karamata's theorem, $\sum_{i=1}^n a_i = O(n
a_n)$. Hence,
\[
{\sup_x} \|\mathcal{P}_k T'_n(x)\|
\le\sum_{j=1}^n 2 \vartheta_{j-k}
\le\sum_{j=1}^{2n} 2 \vartheta_{j}
= o(n^{1-\beta} \ell(n)).
\]
Since $\mathcal{P}_k \cdot= \mathsf{E}( \cdot| \xi_k ) - \mathsf
{E}( \cdot|
\xi_{k-1} )$, $k \in{\mathsf{Z}}$, are orthogonal,
\[
{\sup_x} \|T'_n(x)\|^2 = \sup_x
\Biggl(\sum_{k=-\infty}^{-n} + \sum_{k=1-n}^n \Biggr)
\|\mathcal{P}_k T'_n(x)\|^2 = o(n^{3-2\beta} \ell^2(n)),
\]
where we again applied Karamata's theorem implying
$\sum_{m=n}^\infty m^{-2\beta} \ell^2(m) = O(n^{1-2\beta}
\ell^2(n))$. Similarly, since ${\sup_x} |f'''_{\varepsilon}(x)| <
\infty$, we have $\sup_x \|T''_n(x)\|^2 = o(n^{3-2\beta}
\ell^2(n))$. Since $T'_n(x) = T'_n(0) + \int_0^x T_n''(u) \,d u$,
for all finite $a_0 > 0$,
\[
\mathsf{E} \Bigl[{\sup_{|x|\le a_0}} |T'_n(x)|^2 \Bigr]
= o(n^{3-2\beta} \ell^2(n)).
\]
Hence, (\ref{eq:J12-3}) follows in view of (\ref{eq:aug201}).
\end{pf*}
\begin{pf*}{Proof of Lemma \ref{le2}}
Let $\widetilde{Z}_{k,t} = K (X_{k}/ b-t ) -\mathsf{E}
[K (X_{k}/ b -t ) | \xi_{k-m,k}]$, $Z_{k,t} =
\widetilde{Z}_{k,t} - \mathsf{E}(\widetilde{Z}_{k,t} | \xi_{k-1})$ and
\[
[nb\lambda_K f(b t)]^{1/2}
[M_{n}(t) -\widetilde{M}_{n}(t) - R_{n}(t)]
= \sum_{k=1}^n Z_{k,t}.
\]
We shall approximate $\sum_{k=1}^n Z_{k,t}$ by the skeleton
process $\sum_{k=1}^n Z_{k, t_j}$, $1 \le j \le q_{n}$, where $q_n
= \lfloor n^2/b \rfloor$ and $t_j = j / (b q_n)$. To this end, for
$t\in[t_{j-1}, t_{j}]$, under condition (C4), if $X_{k}/ b -t$
and $X_{k}/ b - t_{j}$ are both in or outside $[-1, 1]$, we have
\[
|K (X_{k}/ b-t )
-K (X_{k}/ b - t_{j} ) | \leq C |t-t_{j}|\leq C n^{-2}.
\]
Otherwise, we have either $|X_{k}/ b -t_{j}-1| \le C n^{-2}$ or
$|X_{k}/ b -t_{j}+1|\leq C n^{-2}$. Let
%
%
\begin{eqnarray}
L_j &=& \sum_{k=1}^{n}I_{kj},\qquad
L^*_j = \sum_{k=1}^{n} \mathsf{E}(I_{kj} | \xi_{k-1}),
\nonumber\\[-8pt]\\[-8pt]
H_j &=& \sum_{k=1}^{n}\mathsf{E}(I_{kj} | \xi_{k-m,k})
\quad\mbox{and}\quad
H^*_j = \sum_{k=1}^{n}\mathsf{E}(I_{kj} | \xi_{k-m,k-1}),\nonumber
\end{eqnarray}
where $I_{k j} = I \{|b^{-1}X_{k}-t_{j}\pm1|\leq Cn^{-2} \}
$. Then
%
%
\begin{equation}\label{p6}
\sup_{t_{j-1}\leq t\leq t_{j}}
\Biggl| \sum_{k=1}^{n}(Z_{k,t}-Z_{k,t_{j}}) \Biggr| \le
{C \over n} + C L_j + C L_j^* + C H_j + C H_j^*.
\end{equation}
Since $f_{X_{n}|\xi_{n-1}}(x)$ is bounded, $\mathsf{E}(I_{kj} |\xi_{k-1}
) \le C n^{-2} b $. Hence, $L_j^* \le C n^{-1} b $ and $D_{k j} =
I_{kj} - \mathsf{E}(I_{kj} |\xi_{k-1} )$ satisfies $\mathsf{E}(D_{k
j}^2 |
\xi_{k-1} )\le C n^{-2} b $. Let $L_\diamond= \max_{1 \leq j\leq
q_{n}} L_j$. Applying the inequality due to \citet{freed75} to
$L_j - L_j^* = \sum_{k=1}^n D_{k j}$, we have
%
%
\begin{eqnarray}\label{p8}\qquad\quad
\mathsf{P} (L_\diamond\geq9 \log n )
&\le&
\mathsf{P} \Bigl({\max_{1 \leq j\leq q_{n}}}
|L_j - L_j^*| \ge8 \log n \Bigr)
+ \mathsf{P} \Bigl(\max_{1 \leq j\leq q_{n}} L_j^*
\geq\log n \Bigr) \nonumber\\[-8pt]\\[-8pt]
&\le& 2 q_n \exp\biggl[ { { (8 \log n)^2}
\over{-2 \times(8 \log n) - 2 C n^{-1} b } }
\biggr] = o(n^{-2}).\nonumber
\end{eqnarray}
Similarly, we have $H_j^* \le C n^{-1} b $, and, for $H_\diamond=
\max_{1 \leq j\leq q_{n}} H_j$, $\mathsf{P}(H_\diamond\geq9 \log n
) =
o(n^{-2})$. Since $\log n = o(\sqrt{n b} / (\log b^{-1})^2)$, by
(\ref{p6}) and (\ref{p8}), it remains to show that
%
%
\begin{equation}\label{p4}
\mathsf{P} \Biggl(\max_{1 \leq j\leq q_{n}} \Biggl|
\sum_{k=1}^{n}Z_{k,t_{j}} \Biggr|
\geq2^{-1}\sqrt{n b} (\log b^{-1})^{-2} \Biggr)=o(1).
\end{equation}

We first consider the case of $X_{n}$ in (\ref{eq:add}). Recall
(\ref{eq:J9-01}) for $\xi^{\star}_{k, n}$. Define
\[
K_{x,t}(\xi_{k-1})=K \biggl(\frac{x+g(\xi_{k-1})}{b}-t \biggr)
\quad\mbox{and}\quad
K^{\Delta}_{x,t}=K_{x,t}(\xi_{k-1})-K_{x,t}(\xi^{\star}_{k-1,m}).
\]
Let $W_k = |g(\xi_{k-1}) - g(\xi^{\star}_{k-1,m})|$. By condition
(C2), $\| W_{k} \|_{p'} = O(m^{-\gamma})$. By Lemma \ref{lem:brk},
we have $ \int_{-\infty}^{\infty} (K^{\Delta}_{x,t} )^{2} \,d x \le
C b \min((W_k / b)^\alpha,1)$. Hence, by Jensen's inequality,
%
%
\begin{eqnarray}\label{p18}\qquad
\mathsf{E}(Z^{2}_{k,t}|\xi_{k-1})&\leq&\int_{-\infty}^{\infty}
\bigl(K_{x,t}(\xi_{k-1})
- \mathsf{E}[K_{x,t}(\xi_{k-1}) | \xi_{k-m,k-1} ] \bigr)^{2}
f_{\varepsilon}(x)\,dx\nonumber\\
&\leq&
\mathsf{E} \biggl[\int_{-\infty}^{\infty} (K^{\Delta}_{x,t}
)^{2}f_{\varepsilon}(x)\,dx \Big| \xi_{k-m,k-1} \biggr]\\
&\leq& C b \mathsf{E} \bigl[\min\bigl((W_k/b)^\alpha,1\bigr) | \xi
_{k-m,k-1} \bigr].\nonumber
\end{eqnarray}
Let $V = \max_{1 \leq j\leq q_n} \sum_{k=1}^n \mathsf{E}(Z^2_{k,
t_j} |
\xi_{k-1})$. Since $\delta_1 / \gamma< \tau< 1 - \delta_1$ and
$m \sim n^\tau$,
%
%
\begin{eqnarray}\label{p5}\quad
\mathsf{P} \biggl( V \geq{ {n b} \over(\log b^{-1})^{6}} \biggr)
&\le& C (\log
b^{-1})^{6}\mathsf{E}\min\bigl((W_k/b)^\alpha,1\bigr)\nonumber\\[-8pt]\\[-8pt]
&\le& C (\log n)^{6} \biggl({ {\Psi_{m,p'}} \over{b}} \biggr)^{\min
(p',\alpha)}
=o(1).\nonumber
\end{eqnarray}
By Freedman's (\citeyear{freed75}) inequality for martingale
differences, we have
\begin{eqnarray*}
&& \mathsf{P} \Biggl(\max_{1 \leq j\leq q_{n}}
\Biggl| \sum_{k=1}^{n}Z_{k,t_{j}} \Biggr|\geq
{ {\sqrt{n b}} \over{2(\log b^{-1})^2}},
V \le{ {n b} \over(\log b^{-1})^{6}} \Biggr)\\
&&\qquad\le2 q_n \exp\biggl[- { {n b (\log b^{-1})^{-4}}
\over{C \sqrt{n b} (\log b^{-1})^{-2}
+ C { {n b} (\log b^{-1})^{-6}} }} \biggr]
= o(1)
\end{eqnarray*}
by condition (C1). So (\ref{p4}) follows from (\ref{p5}).

The proof of (\ref{p4}) for $X_n$ in Theorem \ref{th:22} is
simpler. Let $p_1 = \min(p, 1)$ and $\rho_{1} \in(\rho, 1)$. We
have, by (C2)$'$ and (C3)$'$, that
\begin{eqnarray*}
\sup_{t\in{\mathsf{R}}}\mathsf{E}|Z_{k,t}|
&\leq& C\mathsf{P}(|X_{k}-X^{\star}_{k,m}|\geq\rho^{m}_{1})
+Cb^{-1}\rho^{m}_{1}\\
&&{} +
C\sup_{t\in{\mathsf{R}}}\mathsf{P}(|X_{k}-tb\pm b|\leq\rho^{m}_{1})
\leq C(\rho/\rho^{p_1}_{1})^{m}+Cb^{-1}\rho^{m}_{1}.
\end{eqnarray*}
Hence, using Markov's inequality, (\ref{p4}) follows.
\end{pf*}
\begin{pf*}{Proof of Lemma \ref{le3}}
Let $A=(\log b^{-1})^{-3}=o((\log b^{-1})^{-2})$. Recall the proof
of Lemma \ref{le2} for $t_j$. From the proof of Lemma \ref{le2},
we only need to consider the behavior of $R_{n}(t)$ at grids
$t_j$. Note that $\tau<\tau_{1}$ and
%
%
\begin{equation}\label{p16}\qquad
\sup_{t\in{\mathsf{R}}}\sum_{j=1}^{\iota_n+1}\sum_{k\in
I_{j}}\mathsf{E}\bigl[K^{2}\bigl((X_{k}-t)/b\bigr) |\xi_{k-1}\bigr]\leq
C(n^{1-\tau_{1}+\tau}+n^{\tau_{1}})b\qquad\mbox{a.s.}
\end{equation}
By Freedman's inequality for martingale differences and (\ref{p16}),
\[
\mathsf{P} \Bigl({\max_{0\leq j\leq q_{n}}}|R_{n}(t_{j})|\geq A \Bigr)
\le4 q_n \exp\biggl[ { {A^2 n b}
\over{-2 CA \sqrt{n b} - 2 C(n^{1-\tau_{1}+\tau}+n^{\tau_{1}}) b}}
\biggr] = o(1)
\]
since $n^{-\delta_1} = O(b)$. Hence, (\ref{eq:AugR1}) follows.
\end{pf*}
\begin{pf*}{Proof of Lemma \ref{lem:Aug202}}
From the proof of Lemma \ref{le2}, we only need to show that
\[
{\sup_{0 \le j \le q_{n}}} |M_{n}(t_{j})| = O_\mathsf{P}\bigl(\sqrt
{\log n} \bigr),
\]
which follows from
$\sup_{t\in{\mathsf{R}}}\mathsf{E}[K^{2}((X_{k}-t)/b)|\xi
_{k-1}]\leq Cb$
a.s. and
Freedman's inequality for martingale differences.
\end{pf*}

\subsection{\texorpdfstring{Proof of Lemma \protect\ref{lem:J9-1}}{Proof of Lemma 4.5}}
As in \citet{BR73}, we split the interval $[0,
b^{-1}]$ into alternating big and small intervals $W_{1}, V_{1},
\ldots,W_{N}, V_{N}$, where $W_i = [a_i, a_i + w]$, $V_i = [a_i
+ w, a_{i+1}]$, $a_i = (i-1) (w+v)$, $a_{N+1} =b^{-1}$ and
$N=\lfloor b^{-1}/(w+v)\rfloor$. We will let $v$ be sufficiently
small and $w$ be fixed. We shall first approximate $\Omega^+ :=
\sup_{0\leq t\leq b^{-1}} \widetilde{M}_{n}(t)$ by $\Psi^+ :=
\max_{1\leq k\leq N} \Upsilon_k^+$, where $\Upsilon_k^+ :=
\sup_{t\in W_{k}}\widetilde{M}_{n}(t)$, and then approximate
$\Upsilon_k^+$ via discretization by
%
%
\begin{equation}\label{eq:J134}\qquad
\Xi_k^+ :=\max_{1\le j \le\chi}
\widetilde{M}_{n}(a_{k}+jax^{-2/\alpha})\qquad
\mbox{where }
\chi= \lfloor w x^{2/\alpha}/a \rfloor, a > 0.
\end{equation}
We similarly\vspace*{2pt} define $\Omega^-$, $\Psi^-$, $\Upsilon_k^-$ and
$\Xi_k^-$ by replacing ``$\sup$'' or ``$\max$'' by ``$\inf$'' or
``$\min$,'' respectively. Let $\Omega= \sup_{0\leq t\leq b^{-1}}
|\widetilde{M}_{n}(t)| = \max(\Omega^+, -\Omega^-)$. Define
\begin{eqnarray*}
R_1 &=& \mathsf{P} \Bigl( \max_{1\leq k\leq N}
\sup_{t\in V_{k}}\widetilde{M}_{n}(t)\geq x \Bigr);\qquad
R_2 = \mathsf{P} \Bigl( \min_{1\leq
k\leq N} \inf_{t\in V_{k}}\widetilde{M}_{n}(t)\leq-x \Bigr); \\
R_{3} &=&
\sum_{k=1}^{N}|\mathsf{P} ( \Upsilon_k^+ \ge x
)-\mathsf{P} ( \Xi_k^+ \ge x )|;\\
R_{4} &=& \sum_{k=1}^{N}| \mathsf{P} ( \Upsilon_k^- \le-x
)-\mathsf{P} (\Xi_k^- \le-x )|,
\end{eqnarray*}
where $x = x_z= d_{n} + z / (2\log b^{-1})^{1/2}$. To deal with
$R_{1}, \ldots, R_{4}$, we need the following Lemma \ref{le4}
which will be proved in Section \ref{sec:Aug231}.

Let $(\alpha, C_0)= (1, K_1)$ if $K_1 > 0$ and $(\alpha, C_0)= (2,
K_2)$ if $K_{1}=0$. Let $H_{\alpha}(a)$ and $H_{\alpha}$ be the
Pickands constants [see Theorem A1 and Lemmas A1 and A3 in
\citet{BR73}]. Note that $H_1 = 1$ and $H_2 = 1/\sqrt
\pi$.
\begin{lemma}
\label{le4} Let $t>0$ be such that $\inf\{s^{-\alpha}(1-r(s))\dvtx0
\leq s \leq t \}>0$, where $r(s)$ is defined in Lemma
\ref{lem:brk}. Let $\psi(x)=e^{-x^{2}/2}/(x \sqrt{2\pi})$. Under
conditions of Theorems \ref{th:21} or \ref{th:22}, we have for
$a>0$,
%
%
\begin{eqnarray}\label{eq1}
&&\mathsf{P} \Biggl(\bigcup_{j=1}^{\lfloor t x^{2/\alpha}/a \rfloor}
\{\widetilde{M}_{n}(v+j a x^{-2/\alpha}) \ge x\}
\Biggr)\nonumber\\[-8pt]\\[-8pt]
&&\qquad =x^{2/\alpha}\psi(x)
\frac{H_{\alpha}(a)}{a} C^{1/\alpha}_{0} t
+o(x^{2/\alpha}\psi(x))\nonumber
\end{eqnarray}
uniformly over $0\leq v\leq b^{-1}$. The limit version of
(\ref{eq1}) with $a \to0$ also holds:
%
%
\begin{eqnarray}\label{eq2}
&&\mathsf{P} \Biggl(\bigcup_{0\leq s\leq t}
\{\widetilde{M}_{n}(v+s) \ge x\} \Biggr)\nonumber\\[-8pt]\\[-8pt]
&&\qquad = x^{2/\alpha}\psi(x)H_{\alpha}C^{1/\alpha}_{0}t
+o(x^{2/\alpha}\psi(x)).\nonumber
\end{eqnarray}
The left tail version of (\ref{eq1}) and (\ref{eq2}) also hold
with ``$\ge x$'' replaced by ``$\le-x$.''
\end{lemma}

By Lemma \ref{le4}, elementary calculations show that, for $x = x_z$,
%
%
\begin{equation}\label{eq:J132}
\operatorname{LIM}R_j := \lim_{a\to0}
\limsup_{v \to0} \limsup_{n \to\infty} R_j = 0,
\qquad j=1,\ldots,4.
\end{equation}
Note that $\Omega^+ =\max_{1\leq k\leq N} \sup_{t\in W_{k}\cup
V_{k}}\widetilde{M}_{n}(t)$. By a similar identity for $\Omega^-$,
we have
\[
|\mathsf{P}(\Omega\geq x) -\mathsf{P} ( \{\Psi^+ \geq x\} \cup
\{\Psi^- \le-x\} ) | \le R_{1}+R_{2},
\]
which implies $\mathrm{LIM}|\mathsf{P}(\Omega\geq x) - h(x) | =
0$ for
%
%
\begin{equation}\label{eq:J13hx}
h(x) = \mathsf{P} \Biggl( \bigcup_{k=1}^{N} \{\Xi_k^+ \ge x\} \cup
\bigcup_{k=1}^{N} \{\Xi_k^- \le-x\} \Biggr)
\end{equation}
in view of $|\mathsf{P} ( \{\Psi^+ \geq x\} \cup\{\Psi^- \le
-x\} ) - h(x) | \le R_{3}+R_{4}$. So (\ref{eq3}) follows
from Lemma \ref{lem:J131} below which will be proved in Section
\ref{sec:Aug232}.
\begin{lemma}
\label{lem:J131} Recall (\ref{eq:J132}) for the definition of the
triple limit \textup{LIM}. Under conditions of Theorems \ref{th:21}
or \ref{th:22}, we have $\mathrm{LIM}| h(x_z) - (1-e^{-2e^{-z}})| =
0$ for all $z \in{\mathsf{R}}$.
\end{lemma}

\subsection{\texorpdfstring{Proof of Lemma \protect\ref{le4}}{Proof of Lemma 4.6}}
\label{sec:Aug231} We need the following lemma.
\begin{lemma}[{[Theorems B1 and B2 in \citet{BR73}]}]
\label{lem:brk} Under condition \textup{(C4)}, for $r(s)= \int K(x) K(x+s) \,d x /
\lambda_K$, we have as $s \to0$ that
\[
r(s) = 1 - { {\int(K(x) - K(x+s))^2 \,d x } \over{ 2 \lambda_K}}
= 1-C_{0}|s|^{\alpha}+o(|s|^{\alpha}).
\]
\end{lemma}

Now we prove Lemma \ref{le4}. Assume $C_{0}=1$. The general case
follows from a simple scale transform. Let $s_{j}=j/(\log n)^{6}$,
$1\leq j < t_n$, where $t_n = 1+ \lfloor(\log n)^{6}t \rfloor$,
$s_{t_n} = t$. Write $[s_{j-1}, s_{j}] = \bigcup_{k=1}^{q_{n}}[s_{j,
k-1}, s_{j,k}]$, where $q_{n} = \lfloor(s_{j}-s_{j-1})n^{2}
\rfloor= \lfloor n^2 /(\log n)^{6} \rfloor$ and
$s_{j,k}-s_{j,k-1} = (s_j - s_{j-1}) / q_n$. Define\vspace*{2pt} $\Gamma_j(s) =
\widetilde{M}_n(v+s) - \widetilde{M}_n(v+s_{j-1})$. Using the
arguments in (\ref{p6}) and (\ref{p8}), we have
\[
A_3:= \mathsf{P} \biggl({\max_{1\leq k\leq q_{n}}
\sup_{s_{j,k-1}\leq s\leq s_{j,k}}}
| \Gamma_{j}(s)-\Gamma_{j}(s_{j,k-1})|
> { {(\log n)^{-2}} \over2} \biggr)
\leq{ {C\over{e^{(\log n)^{2}}}}}.
\]
Let $M=2\sqrt{n b}(\log n)^{-4}$. By truncation and Bernstein's inequality,
\begin{eqnarray*}
A_{2}:\!&=& q_{n}\max_{k}
\mathsf{P} \bigl(|\Gamma_{j}(s_{j,k})|>(\log n)^{-2}/2 \bigr)\\
&\leq& q_{n}\max_{k} \biggl[
\exp\biggl( -\frac{Cn b(\log n)^{-4}}{B_{n}} \biggr)
+\exp\biggl(-\frac{C\sqrt{n b}(\log n)^{-2}}{M} \biggr) \biggr]\\
&&{}+q_{n}\mathsf{P} \Biggl( \Biggl|\sum_{l=1}^{\iota_n}(u^{\triangle}_{l}
-\mathsf{E}u^{\triangle}_{l})
\Biggr|\geq\sqrt{n b}(\log n)^{-2}/4 \Biggr),
\end{eqnarray*}
where $u^{\triangle}_{l}=T_l I\{|T_l| \geq\sqrt{n b}(\log
n)^{-4}\}$, $T_l = u_{l}(v+s_{j,k})- u_{l}(v+s_{j-1})$, and
\begin{eqnarray*}
B_{n}&\leq&\sum_{j=1}^{\iota_n}|H_{j}|\mathsf{E} \bigl(
K (X_{1}/b-v-s_{j,k} )-
K (X_{1}/b-v-s_{j-1} ) \bigr)^{2}\\
&\leq& \sum_{j=1}^{\iota_n}|H_{j}| C b |s_{j,k}-s_{j-1}|^\alpha
\le C n b (\log n)^{-6}.
\end{eqnarray*}
Here we applied Lemma \ref{lem:brk}. Since $\tau_{1}<1-\delta_{1}$
and $n^{-\delta_{1}} = O(b)$, for any $Q>2$,
%
%
\begin{equation}\label{p9}
\mathsf{E}|u^{\triangle}_{l}|^2 \leq C(n b)^{-Q/2}(\log n)^{4Q}
n^{\tau_{1} (Q+2)/2}b \leq Cn^{-\tau_{Q}},
\end{equation}
where $\tau_{Q}\to\infty$ as $Q \to\infty$. So $A_{2} \leq
Cn^{-2Q}$ for any $Q>0$, and
\[
A_1 := \mathsf{P} \Bigl({\max_{1\leq j\leq t_{n}}
\sup_{s_{j-1}<s\leq s_{j}}}|\Gamma_{j}(s) |>(\log n)^{-2} \Bigr)
={ {O(t_n)} \over{n^{2Q} }}
\le Cn^{-Q}
\]
for any $Q>0$. Then we have the discretization approximation
\[
\mathsf{P} \Bigl(\sup_{0\leq s\leq t}\widetilde{M}_{n}(v+s)\geq
x \Bigr)
\le\mathsf{P} \Bigl(\max_{1\leq j\leq
t_{n}}\widetilde{M}_{n}(v+s_{j})\geq x-(\log n)^{-2} \Bigr) + A_1.
\]

We now apply the multivariate Gaussian approximation result in
\citet{zait87} to handle $\widetilde{M}_{n}(v)$. To this end,
we introduce
%
%
\begin{eqnarray}\label{p14}
&&\widehat{M}_{n}(t)=\frac{1}{\sqrt{nb\lambda_K f(bt)}}
\sum_{j=1}^{\iota_n}\hat{u}_{j}(t)
\nonumber\\[2pt]\\[-19pt]
&&\eqntext{\mbox{ where }\hat{u}_{j}(t) = u^\diamond_{j}(t) -
\mathsf{E}u^\diamond_{j}(t),\hspace*{80.04pt}}\\
&&\eqntext{u^\diamond_{j}(t)
= u_{j}(t)I\bigl\{|u_{j}(t)|\leq\sqrt{n b}(\log n)^{-20}\bigr\}.}
\end{eqnarray}
As in (\ref{p9}), we have for any large $Q$,
%
%
\begin{equation}\label{p10}
{\sup_{t}\max_{1\leq j\leq\iota_n}}\|\hat{u}_{j}(t)-u_{j}(t)\|
\leq C n^{-Q}.
\end{equation}
By (\ref{p10}) and Theorem 1.1 in \citet{zait87}, we have for
all large $Q$,
%
%
\begin{eqnarray}\label{eq:ttt}\quad
&&\mathsf{P} \Bigl(\max_{1\leq j\leq t_{n}}\widetilde{M}_{n}(v+s_{j})
\geq x-(\log n)^{-2} \Bigr)\nonumber\\
&&\qquad\leq\mathsf{P} \Bigl(\max_{1\leq j\leq t_{n}}
\widehat{M}_{n}(v+s_{j})\geq x-(\log n)^{-2} \Bigr)
+C n^{-Q}\\
&&\qquad\leq\mathsf{P} \Bigl(\max_{1\leq j\leq
t_{n}} Y_{n}(j) \geq x'_n \Bigr)
+ C t^{5/2}_{n}
\exp\biggl(-\frac{C(\log n)^{18}}{t^{5/2}_{n}} \biggr)
+ C n^{-Q},\nonumber
\end{eqnarray}
where $x'_n = x-2(\log n)^{-2}$ and $(Y_{n}(1), \ldots, Y_n(t_n))$
is a centered Gaussian random vector with covariance matrix
%
%
\begin{equation}\label{p10-1}
\widehat{\Sigma}_{n}
=\mathsf{Cov} \bigl(\widehat{M}_{n}(v+s_{1}),\ldots,
\widehat{M}_{n}(v+s_{t_{n}}) \bigr).
\end{equation}
By Lemma \ref{eq:CovJuly} below and Lemma A4 in \citet{BR73}, we have
\begin{eqnarray*}
\mathsf{P} \Bigl(\max_{1\leq j\leq t_{n}}Y_{n}(j)\geq x'_n \Bigr)
&\leq&
\mathsf{P} \Bigl(\max_{1\leq j\leq
t_{n}}\widetilde{Y}_{n}(s_{j})\geq x'_n \Bigr)
+ { {Ct^{2}_{n}(t^{2}_{n}(b+n^{-\varpi}))^{1/2}} \over
{\exp({x'_n}^2/2)}}\\
&\leq& \mathsf{P} \Bigl(\max_{1\leq j\leq
t_{n}}\widetilde{Y}_{n}(s_{j})\geq x'_n \Bigr)+C b^{1+\delta}
\end{eqnarray*}
for some $\delta>0$, where $\widetilde{Y}_{n}(\cdot)$ is a
separable stationary Gaussian process with mean $0$ and covariance
function $r(\cdot)$. By Lemma A3 in \citet{BR73}
and some elementary calculations,
\begin{eqnarray*}
\mathsf{P} \Bigl(\max_{1\leq j\leq t_{n}}\widetilde{Y}_{n}(s_{j})
\geq x'_n \Bigr)
&\leq& \mathsf{P} \Bigl(\sup_{0\leq s\leq t}\widetilde{Y}_{n}(s)
\geq x'_n \Bigr)\\
&=& x^{2/\alpha}\psi(x)H_{\alpha}t
+ o(x^{2/\alpha}\psi(x)).
\end{eqnarray*}
This implies the upper bound in (\ref{eq2}). With the same
argument, for any $a>0$,
\begin{eqnarray*}
&&\mathsf{P} \Bigl(\sup_{0\leq s\leq t}\widetilde{M}_{n}(v+s)\geq
x \Bigr)
\\
&&\qquad\geq \mathsf{P} \Biggl(\bigcup_{j=1}^{[tx^{2/\alpha}/a]}
\{\widetilde{M}_{n}(v+jax^{-2/\alpha})\geq x\} \Biggr)\\
&&\qquad\geq \mathsf{P} \Biggl(\bigcup_{j=1}^{[tx^{2/\alpha}/a]}
\widetilde{Y}_{n}(jax^{-2/\alpha})\geq x
+2(\log n)^{-2} \Biggr)-C b^{1+\delta}\\
&&\qquad\geq\mathsf{P} \Biggl(\bigcup_{j=1}^{[tx^{2/\alpha}/a]}
\widetilde{Y}_{n}(jax^{-2/\alpha})\geq x \Biggr)\\
&&\qquad\quad{}-\sum_{j=1}^{[tx^{2/\alpha}/a]}
\mathsf{P} \bigl(x\leq\widetilde{Y}_{n}(jax^{-2/\alpha})
<x+2(\log n)^{-2} \bigr)
-C b^{1+\delta}\\
&&\qquad=x^{2/\alpha}\psi(x)\frac{H_{\alpha}(a)}{a}t
+o(x^{2/\alpha}\psi(x)).
\end{eqnarray*}
Then the low bound in (\ref{eq2}) is obtained by (A20) in \citet
{BR73}, letting first $n\to\infty$ and then $a\to
0$.

Using a similar and simpler proof, we can prove (\ref{eq1}).
\begin{lemma}
\label{eq:CovJuly} For the covariance matrix $\widehat\Sigma_n$
defined in (\ref{p10-1}), we have
%
%
\begin{equation}\label{p17}\quad
\bigl|\widehat{\Sigma}_{n}-\bigl(r(s_{j}-s_{i})\bigr)_{1\leq i, j\leq t_{n}} \bigr|
\leq C t^{2}_{n}(b+n^{-\varpi}) \qquad\mbox{for some $\varpi>0$.}
\end{equation}
\end{lemma}
\begin{pf}
Let $ \Sigma_{n} = \mathsf{Cov}(\widetilde{M}_{n}(v+s_{1}), \ldots,
\widetilde{M}_{n}(v+s_{t_{n}}))$. By (\ref{p10}), $|\Sigma_{n}
-\widehat{\Sigma}_{n}|\leq Cn^{-Q}$ for any $Q>0$. Note that $\mathsf{E}
(R^{2}_{n}(t) ) \leq Cn^{\tau-\tau_1}$ and $\tau_{1}>\tau$. Then
\[
\bigl|\mathsf{Cov}(\widetilde{M}_{n}(s),\widetilde{M}_{n}(t) )
- \mathsf{Cov}\bigl(\widetilde{M}_{n}(s)+R_{n}(s),
\widetilde{M}_{n}(t)+R_{n}(t) \bigr) \bigr|
\leq C n^{\tau/2 - \tau_1/2}.
\]
By (\ref{p18}), we obtain that $\|\widetilde{M}_{n}(t) + R_{n}(t)
- M_{n}(t)\|^2 \le Cn^{\delta_{1}-\tau\gamma}. $ Thus,
\[
\bigl| \mathsf{Cov} (M_{n}(s),M_{n}(t) )
- \mathsf{Cov}\bigl(\widetilde{M}_{n}(s)+R_{n}(s),
\widetilde{M}_{n}(t)+R_{n}(t)\bigr) \bigr| \leq
Cn^{\delta_{1}/2-\tau\gamma/2}.
\]
Since $K(x)=0$ if $|x| > 1$, for $0\leq s,t\leq b^{-1}$, we have
\[
\bigl| \mathsf{E}[ K ( X_{k}/ b-s )
K ( X_{k}/ b -t ) ]
- b \sqrt{f(b s)f(bt)} r(s-t) \lambda_K \bigr| \leq C b^{2}.
\]
Note that $\mathsf{E}(|K ( X_{k}/ b - t )||\xi_{k-1})\leq
C b$.
Therefore,
\[
|\mathsf{Cov} (M_{n}(s),M_{n}(t) )- r(s-t)| \leq C b.
\]
Combining the above arguments, we prove (\ref{p17}).
\end{pf}

\subsection{\texorpdfstring{Proof of Lemma \protect\ref{lem:J131}}{Proof of Lemma 4.7}}
\label{sec:Aug232} Let $\widehat{M}_{n}(t)$ be defined in
(\ref{p14}) with $20$ therein replaced by $20d$. Also, $d$ may
vary accordingly. Let $x_n = x \pm(\log n)^{-2d}$ and
\begin{eqnarray*}
\mathbf{B}_{k,j}
&=& \{\widetilde{M}_{n}(a_{k}+jax^{-2/\alpha})\geq x\}
\cup\{\widetilde{M}_{n}(a_{k}+jax^{-2/\alpha})\leq-x\},\\
\widehat{\mathbf{B}}^{\pm}_{k,j}
&=&\{\widehat{M}_{n}(a_{k}+jax^{-2/\alpha})
\geq x_n \}\cup\{\widehat{M}_{n}(a_{k}+jax^{-2/\alpha})
\leq-x_n \},\\
\mathbf{D}_{k,j}&=&\{Y_{n}(a_{k}+jax^{-2/\alpha})\geq x\}
\cup\{Y_{n}(a_{k}+jax^{-2/\alpha})\leq-x\},\\
\mathbf{D}^{\pm}_{k,j}&=&\{Y_{n}(a_{k}+jax^{-2/\alpha})\geq x_n\}
\cup\{Y_{n}(a_{k}+jax^{-2/\alpha})\leq-x_n\},\\
\widehat{\mathbf{D}}^{\pm}_{k,j}
&=&\{\widehat{Y}_{n}(a_{k}+jax^{-2/\alpha})\geq x_n \}
\cup\{\widehat{Y}_{n}(a_{k}+jax^{-2/\alpha})\leq-x_n \},
\end{eqnarray*}
where $Y_{n}(\cdot)$ and $\widehat{Y}_{n}(\cdot)$ are centered
Gaussian processes with covariance functions
\begin{eqnarray*}
\mathsf{Cov}(Y_{n}(s_{1}),Y_{n}(s_{2}))
&=& \mathsf{Cov}(\widetilde{M}_{n}(s_{1}),
\widetilde{M}_{n}(s_{2})),\\
\mathsf{Cov}(\widehat{Y}_{n}(s_{1}),
\widehat{Y}_{n}(s_{2})) &=&
\mathsf{Cov}(\widehat{M}_{n}(s_{1}),
\widehat{M}_{n}(s_{2})),
\end{eqnarray*}
respectively. Recall (\ref{eq:J134}) for $\chi$. Let
\[
\mathbf{A}_{k}=\bigcup_{j=1}^\chi\mathbf{B}_{k,j},\qquad
\mathbf{C}_{k}=\bigcup_{j=1}^\chi\mathbf{D}_{k,j},\qquad
\mathbf{C}^{\pm}_{k}=\bigcup_{j=1}^\chi\mathbf{D}^{\pm}_{k,j}
\quad\mbox{and}\quad \widehat{\mathbf{C}}^{\pm}_{k}
=\bigcup_{j=1}^\chi\widehat{\mathbf{D}}^{\pm}_{k,j}.
\]
\begin{lemma}
\label{lem:Aug410} Let $N =\lfloor b^{-1}/(w+v)\rfloor$. Under
the conditions of Theorems \ref{th:21} or~\ref{th:22}, we have for any
fixed integer $l$ satisfying $1\leq l\leq N/2$ that
\[
\Biggl| \mathsf{P} \Biggl(\bigcup_{k=1}^{N}\mathbf{A}_{k} \Biggr)
- \sum_{d=1}^{2l-1}(-1)^{d-1} \biggl(\sum_{1\leq
i_{1}<\cdots<i_{d}\leq N}-\sum_{\mathcal{I}} \biggr)
\mathsf{P} \Biggl(\bigcap_{j=1}^{d}\mathbf{C}_{i_{j}} \Biggr) \Biggr|
\leq
{ {C^{2l}_{1}} \over{(2l)!}} + { {O(1)} \over{\log n}},
\]
where $C_{1}$ does not depend on $l$, and $\mathcal{I}$ is defined in
(\ref{eq7}).
\end{lemma}
\begin{pf} By Bonferroni's inequality, we have
%
%
\begin{eqnarray}\label{eq:bon}
&&\sum_{d=1}^{2l}(-1)^{d-1}
\sum_{1\leq i_{1}<\cdots<i_{d}\leq N}
\mathsf{P} \Biggl(\bigcap_{j=1}^{d}\mathbf{A}_{i_{j}}
\Biggr)\nonumber\\[-8pt]\\[-8pt]
&&\qquad
\leq
\mathsf{P} \Biggl(\bigcup_{k=1}^{N}\mathbf{A}_{k} \Biggr)\leq
\sum_{d=1}^{2l-1}(-1)^{d-1}
\sum_{1\leq i_{1}<\cdots<i_{d}\leq N}
\mathsf{P} \Biggl(\bigcap_{j=1}^{d}\mathbf{A}_{i_{j}} \Biggr).\nonumber
\end{eqnarray}
We now estimate the probability $\mathsf{P}(\bigcap_{j=1}^{d}
\mathbf{A}_{i_j} )$. Recall $W_{k}=[a_{k}, a_{k}+w)$. Let
$q_{j}=i_{j+1}-i_{j}$, $1\leq j\leq d-1$. Define the index set
%
%
\begin{equation}\label{eq7}
\mathcal{I}:= \Bigl\{1\leq i_{1}<\cdots<i_{d}\leq N\dvtx\min_{1\leq
j\leq
d-1}q_{j}\leq\lfloor2w^{-1} + 2 \rfloor\Bigr\}.
\end{equation}
Let $0\leq d_{0}\leq d-2$ and
\[
\mathcal{I}_{d_{0}}= \{1\leq i_{1}<\cdots<i_{d}\leq N\dvtx\mbox{the
number of $j$ such that $q_{j}> \lfloor2w^{-1} + 2 \rfloor$ is
$d_{0}$} \}.
\]
Then we have $\mathcal{I}=\bigcup_{d_{0}=0}^{d-2}\mathcal{I}_{d_{0}}$. We
can see that the number of elements in the sum $ \sum_{\mathcal
{I}_{d_{0}}}\mathsf{P} (\bigcap_{j=1}^{d}\mathbf{A}_{i_{j}}
) $ is
bounded by $CN^{d_{0}+1}=O(b^{-d_{0}-1})$, where $C$ is
independent of $N$. Suppose now $i_{1},\ldots, i_{d}$ are in
$\mathcal{I}_{d_{0}}$. Write
\[
\bigcap_{j=1}^{d}\mathbf{A}_{i_{j}}=\bigcup_{j_{1}=1}^\chi\cdots
\bigcup_{j_{d}=1}^\chi\{\mathbf{B}_{i_{1},j_{1}}\cap\cdots\cap
\mathbf{B}_{i_{d},j_{d}} \}.
\]
Without loss of generality, we assume $q_{1}\leq\lfloor2w^{-1}
+ 2 \rfloor$, $q_2> \lfloor2w^{-1} + 2 \rfloor, \ldots,\break
q_{d_{0}+1}
> \lfloor2w^{-1} + 2 \rfloor$. By (\ref{p10}) and Theorem 1.1 in
\citet{zait87}, we have for all large $Q$,
%
%
\begin{eqnarray}\label{p11}
&&\mathsf{P} (\mathbf{B}_{i_{1},j_{1}}\cap\cdots\cap
\mathbf{B}_{i_{d},j_{d}} )\leq
\mathsf{P} (\widehat{\mathbf{B}}^{-}_{i_{1},j_{1}}\cap\cdots
\cap
\widehat{\mathbf{B}}^{-}_{i_{d},j_{d}} )+Cn^{-Q}\nonumber\\[-8pt]\\[-8pt]
&&\qquad
\leq\mathsf{P} (\widehat{\mathbf{D}}^{-}_{i_{1},j_{1}}
\cap\cdots\cap\widehat{\mathbf{D}}^{-}_{i_{d},j_{d}} )
+C\exp(-(\log b^{-1})^{2})+Cn^{-Q}.\nonumber
\end{eqnarray}
By (\ref{p10}), we have uniformly in $s_{1}$ and $s_{2}$ that, for
any large $Q$,
%
%
\begin{equation}\label{p13}
|\mathsf{Cov}(Y_{n}(s_{1}),Y_{n}(s_{2}))
-\mathsf{Cov}(\widehat{Y}_{n}(s_{1}),\widehat{Y}_{n}(s_{2}))
|\leq Cn^{-Q}.
\end{equation}
Using the argument of (\ref{p17}), there exists $C > 0$ and
$\varpi> 0$, such that for $\nu_n = C(b+n^{-\varpi})$ and any
$1\leq j_{(\cdot)}\leq\chi$, we have
\begin{eqnarray}
\bigl|\mathsf{Cov} \bigl(Y_{n}(a_{i_{l}}+j_{l}ax^{-2/\alpha}),
Y_{n}(a_{i_{k}}+j_{k}ax^{-2/\alpha}) \bigr)\bigr|
&\leq& \nu_n \nonumber\\
\eqntext{\mbox{ for } 3\leq k\leq d_{0}+1, l=1,2;} \\
\bigl|\mathsf{Cov} \bigl(Y_{n}(a_{i_{s}}+j_{s}ax^{-2/\alpha}),
Y_{n}(a_{i_{k}}+j_{k}ax^{-2/\alpha}) \bigr)\bigr|
&\leq& \nu_n\qquad\mbox{for } 3\leq k\neq s\leq d_{0}+1;\nonumber \\
\bigl|\mathsf{Var} \bigl(
Y_{n}(a_{i_{k}}+j_{k}ax^{-2/\alpha}) \bigr)-1\bigr|
&\leq& \nu_n\qquad
\mbox{for } 1\leq k\leq d_{0}+1;\nonumber
\end{eqnarray}
and, letting $\mu= r (a_{i_{2}} -a_{i_{1}} +(j_{2}-j_{1}) a
x^{-2/\alpha})$,
\[
\bigl|\mathsf{Cov} \bigl(Y_{n}(a_{i_{1}}+j_{1}ax^{-2/\alpha}),
Y_{n}(a_{i_{2}}+j_{2}ax^{-2/\alpha}) \bigr)
- \mu\bigr|\le\nu_n.
\]
Note that $|j_{2}-j_{1}|ax^{-2/\alpha} \leq w$ and $a_{i_{2}} -
a_{i_{1}} \geq w+v$ and ${\sup_{x\geq v}}|r(x)| < 1$. Let any $1\leq
j_{(\cdot)}\leq\chi$ and $\mathbf{V}_{n}$ be the covariance
matrix of the Gaussian vector $(\widehat{Y}_{1},\ldots,
\widehat{Y}_{d_{0}+1})$, where $ \widehat{Y}_{k} =
\widehat{Y}_{n}(a_{i_{k}}+j_{k}ax^{-2/\alpha})$, $1\leq k\leq d$.
Using the bounds of the covariances above, we have for some
$\delta>0$ that
%
%
\begin{eqnarray}\label{pt19}\qquad
&& |\mathbf{V}_{n}-\mathbf{V}| \leq Cn^{-\delta}\qquad
\mbox{where } \mathbf{V}
=\pmatrix{
\mathbf{V}_{1} &0\cr
0&\mathbf{I}_{d_{0}-1}
}
\mbox{ and }
\mathbf{V}_{1}=\pmatrix{
1 & \mu\cr
\mu&1}.\hspace*{-12pt}
\end{eqnarray}
By (\ref{pt19}), we have
%
%
\begin{equation}\label{p20}
|\mathbf{V}^{-1}_{n}-\mathbf{V}^{-1} | \leq Cn^{-\delta}\quad \mbox{and}
\quad\bigl|\sqrt{\det(\mathbf{V})}-\sqrt{\det(\mathbf{V}_{n})}\bigr|\leq
Cn^{-\delta}.
\end{equation}
Let $p_{n}(y)$ be the density of
$(\widehat{Y}_{1},\ldots,\widehat{Y}_{d_{0}+1})$, and $p(y)$ be
the density of the Gaussian random vector with covariance matrix
$\mathbf{V}$. By (\ref{p20}), we have
%
%
\begin{eqnarray}\label{p21}\qquad\quad
|p_{n}(y)-p(y)|&\leq& Cn^{-\delta}p(y)+C\exp
(-y\mathbf{V}^{-1}y'/2 )
\bigl|\exp(Cn^{-\delta}|y|^{2} )-1\bigr|\nonumber\\[-8pt]\\[-8pt]
&\leq& C\bigl(n^{-\delta}+n^{-\delta}(\log n)^{2}\bigr)p(y)
+C\exp\bigl(-(\log n)^{2}/C \bigr).\nonumber
\end{eqnarray}
Hereafter, $\delta>0$ may be different in different places. Note
that
\[
|\mu|\leq\sup_{x\geq v}|r(x)|<1.
\]
Then it follows from
Lemma 2 in \citet{berm62} that, for some $\delta>0$, we have
%
%
\begin{eqnarray}\label{p12}
&&\mathsf{P} (\widehat{\mathbf{D}}^{-}_{i_{1},j_{1}}\cap\cdots
\cap\widehat{\mathbf{D}}^{-}_{i_{d},j_{d}} )\nonumber\\
&&\qquad\leq (1+Cn^{-\delta})
\int_{\Xi^{-}}p(y)\,dy+C\exp\bigl(-(\log n)^{2}/C\bigr)\\
&&\qquad\leq C b^{d_{0}+1+\delta},\nonumber
\end{eqnarray}
where $y=(y_{1},\ldots, y_{d_{0}+1})$ and
\[
\Xi^{\pm}=\bigcap_{j=1}^{d_{0}+1}[ \{y_{j}\geq x_n \}
\cup\{y_{j}\le-x_n \} ].
\]
Noting that $\chi^d = O(b^{-\delta/2})$ and by (\ref{p11}) and
(\ref{p12}), we have for some $\delta>0$,
%
%
\begin{equation}\label{p12-1}
\sum_{d_{0}=0}^{d-2}\sum_{\mathcal{I}_{d_{0}}}
\mathsf{P} \Biggl(\bigcap_{j=1}^{d}\mathbf{A}_{i_{j}} \Biggr)\leq
C b^{\delta}.
\end{equation}
We now estimate
%
%
\begin{equation}\label{eq5}
\biggl(\sum_{1\leq i_{1}<\cdots<i_{d}\leq N}-\sum_{\mathcal
{I}} \biggr)\mathsf{P} \Biggl(\bigcap_{j=1}^{d}\mathbf{A}_{i_{j}} \Biggr).
\end{equation}
Suppose that $i_{1},\ldots,i_{d}\notin\mathcal{I}$. Since
$i_{j+1}-i_{j}> \lfloor2/w +2\rfloor$, we have
$a_{i_{j+1}}-a_{i_{j}}\geq(w+v) \lfloor2/w + 2\rfloor>2+w+v$.
Then, for $1\leq s\neq k\leq d$, $1\leq j_{s}, j_{k}\leq\chi$,
\[
\bigl|\mathsf{Cov} \bigl(Y_{n}(a_{i_{s}}+j_{s}ax^{-2/\alpha}),
Y_{n}(a_{i_{k}}+j_{k}ax^{-2/\alpha}) \bigr) \bigr|
\leq C(b +n^{-\varpi})
\]
holds for some $\varpi> 0$. By the bounds of the covariances
above, the covariance matrix $\widetilde{\mathbf{V}}_{n}$ of
$(\widehat{Y}_{1},\ldots, \widehat{Y}_{d})$ when
$i_{1},\ldots,i_{d}\notin\mathcal{I}$ satisfies
%
%
\begin{equation}\label{p19}
|\widetilde{\mathbf{V}}_{n}-\mathbf{I} |
\leq Cn^{-\delta} \qquad\mbox{for some $\delta>0$}.
\end{equation}
For the probability in the sum in (\ref{eq5}), as in (\ref{p11})
and (\ref{p12}), we have for $n$ large,
\begin{eqnarray*}
\mathsf{P} \Biggl(\bigcap_{j=1}^{d}\mathbf{A}_{i_{j}} \Biggr) &\leq&
\sum_{j_{1}=1}^\chi\cdots
\sum_{j_{d}=1}^\chi\mathsf{P} (\mathbf{B}_{i_{1},j_{1}}\cap
\cdots
\cap
\mathbf{B}_{i_{d},j_{d}} )\\
&\leq&
\sum_{j_{1}=1}^\chi\cdots
\sum_{j_{d}=1}^\chi\mathsf{P} (\widehat{\mathbf{D}}^{-}_{i_{1},j_{1}}
\cap\cdots\cap
\widehat{\mathbf{D}}^{-}_{i_{d},j_{d}} )+Cn^{-Q}\\
&\leq&
2^{d}\sum_{j_{1}=1}^\chi\cdots
\sum_{j_{d}=1}^\chi\bigl(x^{-1}\exp(-x^{2}/2) \bigr)^{d}+C
b^{d+\delta}+Cn^{-Q}\\
&\leq&2^{d} \bigl(\chi x^{-1}\exp
(-x^{2}/2) \bigr)^{d}+C
b^{1+\delta} \le C^{d}_{1}b^{d}+C b^{d+\delta}
\end{eqnarray*}
for some $C_{1}>0$ which does not depend on $d$. This together with
(\ref{p12-1}) implies that
%
%
\begin{equation}\label{p19-2}
\sum_{1\leq i_{1}<\cdots<i_{d}\leq N}\mathsf{P} \Biggl(\bigcap
_{j=1}^{d}\mathbf{A}_{i_{j}} \Biggr)\leq C^{d}_{1}/d!+Cb^{\delta}
\end{equation}
for some $C_{1}>0$ which does not depend on $d$. To prove Lemma
\ref{lem:Aug410}, by (\ref{eq:bon}), (\ref{p12-1}) and
(\ref{p19-2}), we only need to show that, for
$i_{1},\ldots,i_{d}\notin\mathcal{I}$,
%
%
\begin{equation}\label{p19-3}
\Biggl|\mathsf{P} \Biggl(\bigcap_{j=1}^{d}\mathbf{A}_{i_{j}} \Biggr)
- \mathsf{P} \Biggl(\bigcap_{j=1}^{d}\mathbf{C}_{i_{j}} \Biggr) \Biggr|
\leq Cb^{d}(\log n)^{-d}.
\end{equation}
By (\ref{p10}) and Theorem 1.1 in \citet{zait87}, as in
(\ref{eq:ttt}), it suffices to show
\[
\Biggl|\mathsf{P} \Biggl(\bigcap_{j=1}^{d}\mathbf{C}_{i_{j}} \Biggr)
- \mathsf{P} \Biggl(\bigcap_{j=1}^{d}\widehat{\mathbf{C}}^{\pm
}_{i_{j}} \Biggr)
\Biggr|\le C b^{d}(\log n)^{-d}.
\]
By (\ref{p13}) and Lemma A4 in \citet{BR73}, using
$ \mathsf{P} (\bigcap_{j=1}^{d}\widehat{\mathbf{C}}^{\pm
}_{i_{j}} )
=1-\mathsf{P} (\bigcup_{j=1}^{d}\widehat{\mathbf{C}}^{\pm
c}_{i_{j}} ) $ and the inclusion--exclusion principle, we have for
any large $Q$,
\[
\Biggl|\mathsf{P} \Biggl(\bigcap_{j=1}^{d}
\widehat{\mathbf{C}}^{\pm}_{i_{j}} \Biggr)
-\mathsf{P} \Biggl(\bigcap_{j=1}^{d}\mathbf{C}^{\pm}_{i_{j}} \Biggr)
\Biggr|\leq C\chi^{2} n^{-2Q}\leq Cn^{-Q}.
\]
So it suffices to show that
%
%
\begin{equation}\label{eq6}
\Biggl|\mathsf{P} \Biggl(\bigcap_{j=1}^{d}\mathbf{C}_{i_{j}} \Biggr)
- \mathsf{P} \Biggl(\bigcap_{j=1}^{d}\mathbf{C}^{\pm}_{i_{j}}
\Biggr) \Biggr|
\leq C b^{d}(\log n)^{-d}.
\end{equation}
By (\ref{p19}) and a similar inequality as (\ref{p21}), we have,
for some $\delta>0$,
\[
|\mathsf{P} ( \mathbf{D}^{\pm}_{i_{1},j_{1}}\cap\cdots
\cap
\mathbf{D}^{\pm}_{i_{d},j_{d}} )- (\mathsf{P}(
\mathbf{D}^{\pm}) )^{d} |\leq C b^{d+\delta},
\]
where $\mathbf{D}^{\pm}= \{\mathbf{N} \ge x_n \} \cup\{
\mathbf{N} \le-x_n \}$ and $\mathbf{N}$ is a standard normal
random variable. It follows that, for some $\delta>0$,
\begin{eqnarray*}
&& \Biggl|\mathsf{P} \Biggl(\bigcap_{j=1}^{d}\mathbf{C}^{-}_{i_{j}} \Biggr)
- \mathsf{P} \Biggl(\bigcap_{j=1}^{d}\mathbf{C}^{+}_{i_{j}} \Biggr)
\Biggr|\\
&&\qquad\leq
\sum_{j_{1}=1}^\chi\cdots\sum_{j_{d}=1}^\chi|
\mathsf{P} ( \mathbf{D}^{-}_{i_{1},j_{1}}\cap\cdots\cap
\mathbf{D}^{-}_{i_{d},j_{d}} )-
\mathsf{P} ( \mathbf{D}^{+}_{i_{1},j_{1}}\cap\cdots\cap
\mathbf{D}^{+}_{i_{d},j_{d}} ) |\\
&&\qquad= \sum_{j_{1}=1}^\chi\cdots
\sum_{j_{d}=1}^\chi| (\mathsf{P}(
\mathbf{D}^{-}) )^{d}
- (\mathsf{P}( \mathbf{D}^{+}) )^{d} |+C
b^{d+\delta}.
\end{eqnarray*}
So (\ref{eq6}) follows from $ \mathsf{P}(\mathbf{D}^{-}) - \mathsf{P}(
\mathbf{D}^{+})\leq C(\log n)^{-2d} b$ and $ \mathsf{P}(\mathbf{D}^{\pm})
\leq C b/\break(\log b^{-1})^{1/\alpha}$.  The lemma is then proved.
\end{pf}

We are ready to prove Lemma \ref{lem:J131}. Let $\{\varepsilon
^{(k)}_{i}\}_{i\in{\mathsf{Z}}}$, $1 \le k \le n$, be i.i.d. copies of
$\{
\varepsilon_{i}\}_{i
\in{\mathsf{Z}}}$, and $\xi^{(k)}_j=(\ldots, \varepsilon^{(k)}_{j-1},
\varepsilon^{(k)}_j)$. Let $X^{(k)}_{j}=G(\xi^{(k)}_j)$. Then
$X^{(k)}_k$, $1\le k\le
n$, are i.i.d. Now define $\mathbf{A}'_k$, $M'_{n}(t)$, $\widetilde
{M}'_{n}(t)$, $N'_{n}(t)$, $R'_{n}(t)$, $R'_1, \ldots, R'_{4}$ by
replacing $X_k$
and $\{\varepsilon_i\}$ by $X^{(k)}_k$ and $\{\varepsilon^{(k)}_i\}$,
respectively, in the above proofs. Repeating the arguments above, we can
obtain that
\[
\Biggl| \mathsf{P} \Biggl(\bigcup_{k=1}^{N}\mathbf{A}'_{k} \Biggr)-
\sum_{d=1}^{2l-1}(-1)^{d-1} \biggl(\sum_{1\leq
i_{1}<\cdots<i_{d}\leq N}-\sum_{\mathcal{I}} \biggr)
\mathsf{P} \Biggl(\bigcap_{j=1}^{d}\mathbf{C}_{i_{j}} \Biggr) \Biggr|
\leq{ {C^{2l}_{1}} \over{(2l)!}} + { {O(1)} \over{\log n}}.
\]
By letting $n\to\infty$ and then $l \to\infty$, we have
\[
\limsup_{n\to\infty} \Biggl|\mathsf{P} \Biggl(\bigcup_{k=1}^{N}\mathbf{A}_{k}
\Biggr)
- \mathsf{P} \Biggl(\bigcup_{k=1}^{N}\mathbf{A}'_{k} \Biggr) \Biggr| =0.
\]
Similarly, (\ref{eq:J132}) holds with $R_j$ therein replaced by
$R_j'$. Hence, as $n\to\infty$,
%
%
\begin{equation}\label{eq:Aug10-1}
\mathrm{LIM} \Biggl|\mathsf{P} \Biggl(\bigcup_{k=1}^{N}\mathbf{A}'_{k}
\Biggr)
- \mathsf{P} \Bigl({\sup_{0\leq t\leq b^{-1}}}|\widetilde
{M}'_{n}(t)|<x \Bigr)
\Biggr|= 0.
\end{equation}
Note that Lemmas \ref{le1}--\ref{le3} also hold for $(X^{(k)}_k)_{ k
\in{\mathsf{Z}}}$, $M'_{n}(t)$, $\widetilde{M}'_{n}(t)$, $N'_{n}(t)$,
$R'_{n}(t)$. By
the theorem in \citet{rosenb76}, the second probability in (\ref
{eq:Aug10-1}) converges to $e^{-2e^{-z}}$. This completes the proof.

\section{\texorpdfstring{Proofs of Proposition \protect\ref{th4},
Theorems \protect\ref{th:scbreg}
and \protect\ref{th:scbreg-2}}{Proofs of Proposition 2.1, Theorems 2.4 and 2.5}}
\label{sec:pfth45}
Without loss of
generality, we assume $l=0$, $u=1$. We first introduce the truncation
\begin{eqnarray*}
\breve{Z}_{k}&=&Z_{k}I\bigl\{|Z_{k}|\leq(\log n)^{{12}/({p-2})}\bigr\}
-\mathsf{E}\bigl(Z_{k}I\bigl\{|Z_{k}|\leq(\log n)^{{12}/({p-2})}\bigr\}\bigr),\\
\widetilde{Z}_{k}&=&Z_{k}I\bigl\{|Z_{k}|>\sqrt{nb}/(\log n)^{4}\bigr\}
-\mathsf{E}\bigl(Z_{k}I\bigl\{|Z_{k}|> \sqrt{nb}/(\log n)^{4}\bigr\}\bigr)
\end{eqnarray*}
and $\widehat{Z}_{k}=Z_{k}-\breve{Z}_{k}$, $1 \le k \le n$.
Correspondingly, define
\begin{eqnarray*}
r_{n}(x)&=&\frac{1}{\sqrt{nb}}\sum_{k=1}^{n}K \biggl(
\frac{X_{k}}{b}-x \biggr)\widehat{Z}_{k}
=:\frac{1}{\sqrt{nb}}\sum_{k=1}^{n}w_{n,k}(x),\\
r_{n,1}(x)&=&\frac{1}{\sqrt{nb}}\sum_{k=1}^{n}K \biggl(
\frac{X_{k}}{b}-x \biggr)\widetilde{Z}_{k}
=:\frac{1}{\sqrt{nb}}\sum_{k=1}^{n}w_{n,k1}(x),\\
r_{n,2}(x)&=&r_{n}(x)-r_{n,1}(x)
=:\frac{1}{\sqrt{nb}}\sum_{k=1}^{n}w_{n,k2}(x).
\end{eqnarray*}
\begin{lemma}
\label{lem:r} Under the conditions of Proposition \ref{th4}, we have
\[
\mathsf{P} \Bigl({\sup_{0\leq x\leq b^{-1}}}|r_{n}(x)|
\geq3(\log n)^{-2} \Bigr)=o(1).
\]
\end{lemma}
\begin{pf}
Since $b\geq Cn^{-\delta_{1}}$ and $\mathsf{E}|Z_{1}|^{p}<\infty$,
$p>2/(1-\delta_{1})$, for $n$ large, we have
%
%
\begin{eqnarray}\label{eq:Aug6-2}
{\mathsf{E}\sup_{0\leq x\leq b^{-1}}}|r_{n,1}(x)|
&\leq& C n(n b)^{-p/2}(\log n)^{4p-4}\nonumber\\[-8pt]\\[-8pt]
&\leq& Cn^{1-p(1-\delta_{1})/2}(\log n)^{4p-4}
\leq(\log n)^{-3} .\nonumber
\end{eqnarray}
We now deal with $r_{n,2}$. Let $q_n = \lfloor n^2/b \rfloor$,
$t_j = j / (b q_n)$, $j = 0, \ldots, q_n$. As in (\ref{p6}), we have
%
%
\begin{equation}\label{p6-6}\quad
{\max_{0\leq j\leq q_{n}}\sup_{t_{j}\leq t\leq t_{j+1}}}
| r_{n,2}(t)-r_{n,2}(t_{j}) | \leq{ C \over{n (\log n)^4} }
+ C { {\max_{0\leq j\leq q_{n}} L_j} \over{(\log n)^4}}.
\end{equation}
By (\ref{p8}), (\ref{eq:Aug6-2}), (\ref{p6-6}) and since
$r_{n,2}(x) + r_{n,1}(x) = r_{n}(x)$, it suffices to show
%
%
\begin{equation}\label{eq8}
\mathsf{P} \Bigl({\max_{0\leq j\leq q_{n}}}|r_{n,2}(t_{j})|
\geq2(\log n)^{-2} \Bigr)=o(1).
\end{equation}
Note that $\mathsf{E}(\widehat{Z}^{2}_{k}) \leq C(\log n)^{-12}$. By
(C3) [or (C3)$'$], we have
%
%
\begin{equation}\label{eq9}
\max_{0\leq j\leq q_{n}}\sum_{k=1}^{n}
\mathsf{E} [w^{2}_{n,k2}(t_{j}) | \widetilde{\xi}_{k-2} ]
\leq Cnb(\log n)^{-6}.
\end{equation}
Thus, (\ref{eq8}) follows from (\ref{eq9}) and applying Freedman's
inequality to martingale differences $\{w_{n,k2}(x),
k=1,3,\ldots\}$ and $\{w_{n,k2}(x), k=2,4,\ldots\}$.
\end{pf}
\begin{pf*}{Proof of Proposition \ref{th4}}
Let $m = \lfloor n^{\tau} \rfloor$, where
$\delta_{1}/\gamma<\tau<1-\delta_{1}$, and
\[
Z_{k}(t)=\breve{Z}_{k} \biggl\{K \biggl(\frac{X_{k}}{b}-t \biggr)
-\mathsf{E} \biggl[K \biggl(\frac{X_{k}}{b}-t \biggr)
\Big|\xi_{k-m,k} \biggr] \biggr\},\qquad
1\leq k\leq n.
\]
Note that $\{Z_{1}(t), Z_{3}(t), \ldots\}$ and $\{Z_{2}(t),
Z_{4}(t),\ldots\}$ are two sequences of martingale differences. As
in the proof of Lemma \ref{le2}, we can show that
%
%
\begin{eqnarray}\label{eq10-1}
\mathsf{P} \Biggl(\sup_{0\leq t\leq b^{-1}} \Biggl|
\sum_{k=1}^{n/2}Z_{2k-1}(t) \Biggr|
\geq\sqrt{nb}(\log n)^{-2} \Biggr)&=&o(1),\nonumber\\[-8pt]\\[-8pt]
\mathsf{P} \Biggl(\sup_{0\leq t\leq b^{-1}}
\Biggl|\sum_{k=1}^{n/2}Z_{2k}(t) \Biggr|
\geq\sqrt{nb}(\log n)^{-2} \Biggr)&=&o(1).\nonumber
\end{eqnarray}
Set
\[
\widetilde{N}_{n}(t)
=\frac{1}{\sqrt{nb\lambda_K f(bt)}}\sum_{k=1}^{n}\mathsf{E} \biggl[
K \biggl(\frac{X_{k}}{b}-t \biggr) \Big|\xi_{k-m,k-1} \biggr]\breve{Z}_{k}.
\]
Since $\sup_{t}\mathsf{E}(\{\breve{Z}_{k}\mathsf{E}[ K (X_{k}
/ b - t
) | \xi_{k-m,k-1}]\}^{2}|\widetilde{\xi}_{k-1} ) \leq Cb^{2}$, we
have by
Freedman's inequality for martingale differences,
\[
\mathsf{P} \Bigl({\max_{0\leq j\leq q_{n}}}
|\widetilde{N}_{n}(t_{j})|\geq(\log n)^{-2} \Bigr)=o(1),
\]
which, together with the discretization approximation as in
(\ref{p6}), yields that
%
%
\begin{equation}\label{eq10-2}
\mathsf{P} \Bigl({\sup_{0\leq t\leq b^{-1}}}|\widetilde{N}_{n}(t)|
\geq2(\log n)^{-2} \Bigr)=o(1).
\end{equation}
Set $\breve{\sigma}^{2}_n=\mathsf{E}\breve{Z}^{2}_{n}$ and
\begin{eqnarray*}
\widetilde{M}_{n}(t)&=&\frac{1}{\sqrt{nb\lambda_K f(bt)}}\\
&&{}\times\sum_{k=1}^{n} \biggl\{\mathsf{E} \biggl[K \biggl(\frac
{X_{k}}{b}-t \biggr)
\Big|\xi_{k-m,k} \biggr]-\mathsf{E} \biggl[ K \biggl(\frac{X_{k}}{b}
-t \biggr) \Big|\xi_{k-m,k-1} \biggr]
\biggr\}
\frac{\breve{Z}_{k}}{\breve{\sigma}_n}.
\end{eqnarray*}
Following the argument of Lemma \ref{lem:J9-1} and replacing the
truncation levels $(\log n)^{-20}$ and $(\log n)^{-20d}$ in
(\ref{p14}) and the proof of Lemma \ref{lem:J131} with\break $(\log
n)^{- 20p / (p-2)}$ and $(\log n)^{- 20pd/(p-2)}$, respectively, we
can get
%
%
\begin{equation}\label{eq10-3}
\mathsf{P} \Bigl((2\log b^{-1})^{1/2}
\Bigl({\sup_{0\leq t\leq b^{-1}}}
|\widetilde{M}_{n}(t)| -d_{n} \Bigr)
\leq z \Bigr) \rightarrow e^{-2e^{-z}}.
\end{equation}
Note that $|1-\breve{\sigma}^2_n/\sigma^2|=O((\log n)^{-12})$. The
proposition follows from Lemma \ref{lem:r} and
(\ref{eq10-1})--(\ref{eq10-3}).
\end{pf*}
\begin{pf*}{Proof of Theorem \ref{th:scbreg}}
Write $(\mu_{n}(x) -\mu(x)) f_{n}(x) = R^{r}_{n}(x) +
M^{r}_{n1}(x)$, where
\begin{eqnarray*}
R^{r}_{n}(x)&=&\frac{1}{nb}\sum_{k=1}^{n}
K \biggl(\frac{X_{k}-x}{b} \biggr)
\bigl( \mu(X_{k})-\mu(x) \bigr), \\
M^{r}_{n1}(x)&=& \frac{1}{nb}\sum_{k=1}^{n}
K \biggl(\frac{X_{k}-x}{b} \biggr) \sigma(X_{k})\eta_{k}.
\end{eqnarray*}
Then Theorem \ref{th:scbreg} follows from Lemmas \ref{lem:Aug202},
\ref{lem:s} and \ref{lem:j} and Proposition \ref{th4}.
\end{pf*}
\begin{lemma}\label{lem:s}
Under the conditions of Theorem \ref{th:scbreg}, we have
\[
{\sup_{0\leq x\leq1}} |R^{r}_{n}(x)-b^{2}\psi_K\rho_{\mu}(x)|
=O_{\mathsf{P}}(\tau_{n})\qquad\mbox{where }
\tau_{n}=\sqrt{{{b\log n} \over n}}+b^{4}+
{ {\mathcal{Z}^{1/2}_{n}b} \over n}.
\]
\end{lemma}
\begin{pf}
Set $\gamma_{k}(x)=K((X_{k}-x)/b)(\mu(X_{k}) - \mu(x))$. Let $q_n
= \lfloor n^2/b \rfloor$, $t_j = j / q_n$, $j = 0, \ldots, q_n$.
Since $\mu(\cdot) \in\mathcal{C}^{4} (T^{\epsilon})$, $\max_{0\le
j \le q_{n}} \mathsf{E}[\gamma^{2}_{k}(t_{j})|\xi_{k-1}] \le C
b^{3}$. By
Freedman's inequality for martingale differences, we have
\[
\max_{0\leq j\leq q_{n}}
\Biggl|\sum_{k=1}^{n}\bigl(\gamma_{k}(t_{j}) -
\mathsf{E}[\gamma_{k}(t_{j}) | \xi_{k-1}]\bigr) \Biggr|
= O_{\mathsf{P}}\bigl(\sqrt{nb^{3}\log n}\bigr),
\]
where we used the condition $0 < \delta_{1} < 1/3$. Recall that
$K(x)$ and $m(x)$ are Lipschitz continuous in $[-1,1]$. Using the
discretization approximation as in (\ref{p6}) and the argument in
(\ref{p8}), it can be seen that
\[
\sup_{0\leq x\leq1} \Biggl|\sum_{k=1}^{n}\bigl(\gamma_{k}(x)
-\mathsf{E}[\gamma_{k}(x)|\xi_{k-1}]\bigr) \Biggr|
= O_{\mathsf{P}}\bigl(\sqrt{nb^{3}\log n}\bigr).
\]
The rest of the proof is the same as that of Lemma 2(ii) in \citet{ZW08}.
\end{pf}
\begin{lemma}\label{lem:j}
Under the conditions of Theorem \ref{th:scbreg}, we have
\[
\sup_{0\leq x\leq1} \Biggl|M^{r}_{n1}(x)-\frac{1}{nb}
\sum_{k=1}^{n}K \biggl(\frac{X_{k}-x}{b} \biggr)\sigma(x)\eta_{k}
\Biggr| =O_{\mathsf{P}}\Biggl(\sqrt{{{b\log n} \over n}}\Biggr).
\]
\end{lemma}
\begin{pf} Let
\begin{eqnarray*}
\widetilde{\eta}_{k}
&=&\eta_{k}I\bigl\{|\eta_{k}|\geq\sqrt{nb}/(\log n)^{4}\bigr\}
-\mathsf{E}\bigl(\eta_{k}I\bigl\{|\eta_{k}|\geq\sqrt{nb}/(\log n)^{4}\bigr\}\bigr),\\
\widetilde{w}_{nk}(x)&=&K \biggl(\frac{X_{k}-x}{b} \biggr)
\bigl(\sigma(X_{k})-\sigma(x)\bigr)\widetilde{\eta}_{k},\\
\widehat{w}_{nk}(x)&=&K \biggl(\frac{X_{k}-x}{b} \biggr)
\bigl(\sigma(X_{k})-\sigma(x)\bigr)\widehat{\eta}_{k},\qquad
\widehat{\eta}_{k}=\eta_{k}-\widetilde{\eta}_{k}.
\end{eqnarray*}
Note that ${\sup_{x\in T^\epsilon}} |K((X_{k}-x)/b)
(\sigma(X_{k})-\sigma(x))| \le C b$. Then
\[
\mathsf{E}\sup_{x\in{\mathsf{R}}} \Biggl|\frac{1}{nb}
\sum_{k=1}^{n}\widetilde{w}_{nk}(x) \Biggr|
=O\Biggl(\sqrt{{{b} \over n(\log n)^{4}}}\Biggr).
\]
Since $\sup_{x \in{\mathsf{R}}} \mathsf{E}[\widehat{w}^{2}_{nk}(x) |
\widetilde{\xi}_{k-2} ] \leq C b^{3}$, we have
\[
\sup_{x\in{\mathsf{R}}}\sum_{k=1}^{n}
\mathsf{E} [\widehat{w}^{2}_{nk}(x) |\widetilde{\xi
}_{k-2} ]
\leq Cnb^{3}.
\]
Using the arguments for (\ref{p6-6}) and (\ref{eq8}), we can show that
\[
\sup_{0\leq x\leq1} \Biggl|\frac{1}{nb}
\sum_{k=1}^{n}\widehat{w}_{nk}(x) \Biggr|
=O_{\mathsf{P}}\Biggl(\sqrt{{{b\log n} \over n}}\Biggr).
\]
The lemma is proved.
\end{pf}
\begin{pf*}{Proof of Theorem \ref{th:scbreg-2}}
Write
%
%
\begin{eqnarray}\label{eq10-4}
\sigma^{2}_{n}(x)&=&\frac{1}{n h f_{n1}(x)}
\sum_{k=1}^{n}K \biggl(\frac{X_{k}-x}{h} \biggr)
[\sigma(X_{k})\eta_{k}]^{2}\nonumber\\
&&{}+\frac{2}{n h f_{n1}(x)}
\sum_{k=1}^{n}K \biggl(\frac{X_{k}-x}{h} \biggr)
[\mu(X_{k})-\mu_{n}(X_{k})]\sigma(X_{k})\eta_{k}\nonumber\\[-8pt]\\[-8pt]
&&{}+\frac{1}{n h f_{n1}(x)}
\sum_{k=1}^{n}K \biggl(\frac{X_{k}-x}{h} \biggr)[\mu(X_{k})-\mu
_{n}(X_{k})]^{2}\nonumber\\
&=&\!:\sigma^{2}_{n1}(x)+c_{n2}(x)+\sigma^{2}_{n3}(x).\nonumber
\end{eqnarray}
We have
%
%
\begin{eqnarray}\label{eq10-5}
{\sup_{0\leq x\leq1}}|\sigma^{2}_{n3}(x)|
&=&O_{\mathsf{P}}\biggl(\frac{\log n}{nb}+b^{4}\biggr)\nonumber\\
&&{}\times\sup_{0\leq x\leq1}\frac{1}{nh}
\sum_{k=1}^{n} \biggl|K \biggl(\frac{X_{k}-x}{h} \biggr)
\biggr|\\
&=&O_{\mathsf{P}}\biggl(\frac{\log n}{nb}+b^{4}\biggr).\nonumber
\end{eqnarray}
Using a similar argument as in Zhao and Wu [(\citeyear{ZW08}), page 1875] we have
%
%
\begin{equation}\label{eq10-6}
{\sup_{0\leq x\leq1}}|c_{n2}(x)|=O_{\mathsf{P}}\biggl(\frac{1}{n b^{5/2}}\biggr).
\end{equation}
For $\sigma^{2}_{n1}(x)$,
%
%
\begin{eqnarray}\label{eq10-7}
&&\bigl(\sigma^{2}_{n1}(x)-\sigma^{2}(x)\bigr)f_{n1}(x)
\nonumber\\
&&\qquad=\frac{1}{nh}\sum_{k=1}^{n}K \biggl(\frac{X_{k}-x}{h}
\biggr)\sigma^{2}(x)
(\eta^{2}_{k}-1)\nonumber\\
&&\qquad\quad{}+
\frac{1}{nh}\sum_{k=1}^{n}K \biggl(\frac{X_{k}-x}{h} \biggr)
\bigl(\sigma^{2}(X_{k})-\sigma^{2}(x)\bigr)(\eta^{2}_{k}-1)\\
&&\qquad\quad{}+
\frac{1}{nh}\sum_{k=1}^{n}K \biggl(\frac{X_{k}-x}{h} \biggr)
\bigl(\sigma^{2}(X_{k})-\sigma^{2}(x)\bigr)\nonumber\\
&&\qquad=:M^{r}_{n2}(x)+R^{r}_{n2}(x)+R^{r}_{n3}(x).\nonumber
\end{eqnarray}
As in the proof of Lemma \ref{lem:j}, we get
%
%
\begin{equation}\label{eq10-8}
{\sup_{0\leq x\leq1}}|R^{r}_{n2}(x)|=O_{\mathsf{P}}\Biggl(\sqrt{{{b\log n}
\over n}}\Biggr).
\end{equation}
Also, for $R^{r}_{n2}(x)$, we have similarly as in Lemma \ref{lem:s} that
%
%
\begin{equation}\label{eq10-9}
{\sup_{0\leq x\leq1}} |R^{r}_{n2}(x)-h^{2}\psi_K\rho_{\sigma}(x)|
=O_{\mathsf{P}}(\tau_{n}).
\end{equation}
Theorem \ref{th:scbreg-2} now follows from Lemma \ref{lem:Aug202},
Proposition \ref{th4} and (\ref{eq10-4})--(\ref{eq10-9}).
\end{pf*}

\section*{Acknowledgments}
We are grateful to two referees and an
Associate Editor for their many helpful comments.

\printaddresses


\begin{thebibliography}{99}

\bibitem[\protect\citeauthoryear{A\"{i}t-Sahalia}{1996a}]{Ait1}
\textsc{A\"{i}t-Sahalia, Y.} (1996a). Nonparametric pricing of interest rate
derivative securities. \textit{Econometrica} \textbf{64} 527--560.

\bibitem[\protect\citeauthoryear{A\"{i}t-Sahalia}{1996b}]{Ait2}
\textsc{A\"{i}t-Sahalia, Y.} (1996b). Testing continuous-time models of the
spot interest rate. \textit{Rev. Finan. Stud.} \textbf{9} 385--426.

\bibitem[\protect\citeauthoryear{Berman}{1962}]{berm62}
\textsc{Berman, S.} (1962). A law of large numbers for the maximum of a
stationary Gaussian sequence. \textit{Ann. Math. Statist.} \textbf{33}
93--97.
\MR{0133856}

\bibitem[\protect\citeauthoryear{Bickel and Rosenblatt}{1973}]{BR73}
\textsc{Bickel, P. J.} and \textsc{Rosenblatt, M.} (1973).
On some global measures of the deviations of density
function estimates. \textit{Ann. Statist.}
\textbf{1} 1071--1095.
\MR{0348906}

\bibitem[\protect\citeauthoryear{Black and Scholes}{1973}]{BS73}
\textsc{Black, F.} and \textsc{Scholes, M.} (1973). The pricing of
options and
corporate liabilities. \textit{Journal of Political Economy} \textbf
{81} 637--654.

\bibitem[\protect\citeauthoryear{Bosq}{1996}]{bosq96}
\textsc{Bosq, D.} (1996). \textit{Nonparametric Statistics for Stochastic
Processes. Estimation and Prediction.}
\textit{Lecture Notes in Statistics}
\textbf{110}. Springer, New York.
\MR{1441072}

\bibitem[\protect\citeauthoryear{Brillinger}{1969}]{brill69}
\textsc{Brillinger, D. R.} (1969). An asymptotic representation of the
sample distribution function. \textit{Bull. Amer. Math. Soc.} \textbf
{75} 545--547.
\MR{0243659}

\bibitem[\protect\citeauthoryear{B{\"{u}}hlmann}{1998}]{buhl98}
\textsc{B{\"{u}}hlmann, P.} (1998). Sieve bootstrap for smoothing in
nonstationary time series. \textit{Ann. Statist.} \textbf{26} 48--83.
\MR{1611804}

\bibitem[\protect\citeauthoryear{Chan et al.}{1992}]{chanetal92}
\textsc{Chan, K. C., Karolyi, A. G., Longstaff, F. A.} and \textsc
{Sanders, A.
B.} (1992). An empirical comparison of alternative models of the
short-term interest rate. \textit{J. Finance} \textbf{47} 1209--1227.

\bibitem[\protect\citeauthoryear{Chapman and Pearson}{2000}]{ChP00}
\textsc{Chapman, D. A.} and \textsc{Pearson, N. D.} (2000). Is the short rate drift
actually nonlinear? \textit{J. Finance} \textbf{55} 355--388.

\bibitem[\protect\citeauthoryear{Courtadon}{1982}]{court82}
\textsc{Courtadon, G.} (1982). The pricing of options on default-free
bonds. \textit{J. Finan. Quant. Anal.} \textbf{17} 75--100.

\bibitem[\protect\citeauthoryear{Cox, Ingersoll and Ross}{1985}]{CIR85}
\textsc{Cox, J. C., Ingersoll, J. E.} and \textsc{Ross, S. A.} (1985).
A theory of
the term structure of interest rates. \textit{Econometrica} \textbf{53}
385--403.
\MR{0785475}


\bibitem[\protect\citeauthoryear{Cummins, Filloon and Nychka}{2001}]{CFN01}
\textsc{Cummins, D. J., Filloon, T. G.} and \textsc{Nychka, D.} (2001).
Confidence
intervals for nonparametric curve estimates: Toward more uniform
pointwise coverage. \textit{J. Amer. Statist.
Assoc.} \textbf{96} 233--246.
\MR{1952734}

\bibitem[\protect\citeauthoryear{Doukhan and Louhichi}{1999}]{DL99}
\textsc{Doukhan, P.} and \textsc{Louhichi, S.} (1999). A new weak dependence
condition and applications to moment inequalities. \textit{Stochastic
Process. Appl.} \textbf{84} 313--342.
\MR{1719345}

\bibitem[\protect\citeauthoryear{Doukhan, Madre and Rosenbaum}{2007}]{DMR07}
\textsc{Doukhan, P., Madre, H.} and \textsc{Rosenbaum, M.} (2007). Weak
dependence
for infinite ARCH-type bilinear models. \textit{Statistics} \textbf{41}
31--45.
\MR{2303967}

\bibitem[\protect\citeauthoryear{Doukhan and Portal}{1987}]{DP87}
\textsc{Doukhan, P.} and \textsc{Portal, F.} (1987). Principe
d'invariance faible
pour la fonction de r\'{e}partition empirique dans un cadre
multidimensionnel et m\'{e}langeant. \textit{Probab. Math. Statist.}
\textbf{8} {117--132}.
\MR{0928125}

\bibitem[\protect\citeauthoryear{D\"{u}mbgen}{2003}]{dumb03}
\textsc{D\"{u}mbgen, L.} (2003). Optimal confidence bands for
shape-restricted curves. \textit{Bernoulli} \textbf{9} 423--449.
\MR{1997491}

\bibitem[\protect\citeauthoryear{Erd{\"{o}}s}{1939}]{erdos39}
\textsc{Erd{\"{o}}s, P.} (1939). On a family of symmetric Bernoulli
convolutions. \textit{Amer. J. Math.} \textbf{61} 974--976.
\MR{0000311}


\bibitem[\protect\citeauthoryear{Fan and Yao}{1998}]{fan98}
\textsc{Fan, J.} and \textsc{Yao, Q.} (1998). Efficient estimation of
conditional
variance functions in stochastic regression.
\textit{Biometrika} \textbf{85} 645--660.
\MR{1665822}

\bibitem[\protect\citeauthoryear{Fan and Yao}{2003}]{fan03}
\textsc{Fan, J.} and \textsc{Yao, Q.} (2003). \textit{Nonlinear Time Series.
Nonparametric and Parametric Methods.} Springer, New York.
\MR{1964455}

\bibitem[\protect\citeauthoryear{Fan and Zhang}{2003}]{FZ03}
\textsc{Fan, J.} and \textsc{Zhang, C.} (2003). A re-examination of diffusion
estimators with applications to financial model validation.
\textit{J. Amer. Statist. Assoc.} \textbf{98} 118--134.
\MR{1965679}

\bibitem[\protect\citeauthoryear{Freedman}{1975}]{freed75}
\textsc{Freedman, D. A.} (1975). On tail probabilities for martingales.
\textit{Ann. Probab.} \textbf{3} 100--118.
\MR{0380971}

\bibitem[\protect\citeauthoryear{Granger and Joyeux}{1980}]{GJ80}
\textsc{Granger, C. W. J.} and \textsc{Joyeux, R.} (1980). An
introduction to
long-memory time series models and fractional differencing.
\textit{J. Time Ser. Anal.} \textbf{1} 15--29.
\MR{0605572}

\bibitem[\protect\citeauthoryear{Gy\"{o}rfi et al.}{1989}]{gyorfietal89}
\textsc{Gy\"{o}rfi, L., H\"{a}rdle, W., Sarda, P.} and \textsc{Vieu,
P.} (1989).
\textit{Nonparametric Curve Estimation From Time Series.} Springer,
Berlin.
\MR{1027837}

\bibitem[\protect\citeauthoryear{H\"{a}rdle and Marron}{1991}]{HM91}
\textsc{H\"{a}rdle, W.} and \textsc{Marron, J. S.} (1991). Bootstrap
simultaneous
error bars for nonparametric regression. \textit{Ann. Statist.} \textbf
{19} 778--796.
\MR{1105844}

\bibitem[\protect\citeauthoryear{Hall and Titterington}{1988}]{HT88}
\textsc{Hall, P.} and \textsc{Titterington, D. M.} (1988). On
confidence bands in
nonparametric density estimation and regression.
\textit{J. Multivariate Anal.} \textbf{27} 228--254.
\MR{0971184}

\bibitem[\protect\citeauthoryear{Ho and Hsing}{1996}]{HH96}
\textsc{Ho, H. C.} and \textsc{Hsing, T.} (1996). On the asymptotic
expansion of the
empirical process of long-memory moving averages. \textit{Ann.
Statist.} \textbf{24} 992--1024.
\MR{1401834}

\bibitem[\protect\citeauthoryear{Hosking}{1981}]{hosking81}
\textsc{Hosking, J. R. M.} (1981).
Fractional differencing. \textit{Biometrika} \textbf{68} 165--176.
\MR{0614953}

\bibitem[\protect\citeauthoryear{Johnston}{1982}]{john82}
\textsc{Johnston, G. J.} (1982). Probabilities of maximal deviations for
nonparametric regression function estimates.
\textit{J. Multivariate Anal.} \textbf{12} 402--414.
\MR{0666014}

\bibitem[\protect\citeauthoryear{Knafl, Sacks and Ylvisaker}{1985}]{KSY85}
\textsc{Knafl, G., Sacks, J.} and \textsc{Ylvisaker, D.} (1985).
Confidence bands
for regression functions. \textit{J.~Amer. Statist.
Assoc.} \textbf{80} 683--691.
\MR{0803261}

\bibitem[\protect\citeauthoryear{Koml\'{o}s, Major and Tusn\'
{a}dy}{1975}]{KMT75}
\textsc{Koml\'{o}s, J., Major, P.} and \textsc{Tusn\'{a}dy, G.} (1975).
An approximation of
partial sums of independent $\mathrm{RV}$'s and the sample $\mathrm{DF}$.
I. \textit{Z. Wahrsch. Verw. Gebiete} \textbf{32} 111--131.
\MR{0375412}

\bibitem[\protect\citeauthoryear{Koml\'{o}s, Major and Tusn\'
{a}dy}{1976}]{KMT76}
\textsc{Koml\'os, J., Major, P.} and \textsc{Tusn\'ady, G.} (1976).
An approximation of partial sums of independent RV's and the sample DF. II.
\textit{Z. Wahrsch. Verw. Gebiete} \textbf{34} 33--58.
\MR{0402883}

\bibitem[\protect\citeauthoryear{Neumann}{1998}]{neum98}
\textsc{Neumann, M. H.} (1998). Strong approximation of density estimators
from weakly dependent observations by density estimators from
independent observations. \textit{Ann. Statist.} \textbf{26} 2014--2048.
\MR{1673288}

\bibitem[\protect\citeauthoryear{Robinson}{1983}]{robin83}
\textsc{Robinson, P. M.} (1983). Nonparametric estimators for time series.
\textit{J. Time Ser. Anal.} \textbf{4} 185--207.
\MR{0732897}

\bibitem[\protect\citeauthoryear{Rosenblatt}{1976}]{rosenb76}
\textsc{Rosenblatt, M.} (1976). On the maximal deviation of $k$-dimensional
density estimates. \textit{Ann. Probab.} \textbf{4} 1009--1015.
\MR{0428580}

\bibitem[\protect\citeauthoryear{Stanton}{1997}]{stant97}
\textsc{Stanton, R.} (1997). A nonparametric model of term structure
dynamics and the market price of interest rate risk.
\textit{J. Finance} \textbf{52} 1973--2002.

\bibitem[\protect\citeauthoryear{Shao and Wu}{2007}]{SW07}
\textsc{Shao, X.} and \textsc{Wu, W. B.} (2007). Asymptotic spectral
theory for
nonlinear time series. \textit{Ann. Statist.} \textbf{35} 1773--1801.
\MR{2351105}

\bibitem[\protect\citeauthoryear{Sun and Loader}{1994}]{SL94}
\textsc{Sun, J.} and \textsc{Loader, C. R.} (1994). Simultaneous
confidence bands
for linear regression and smoothing. \textit{Ann. Statist.} \textbf{22}
1328--1345.
\MR{1311978}

\bibitem[\protect\citeauthoryear{Tj{\o}stheim}{1994}]{tjos94}
\textsc{Tj{\o}stheim, D.} (1994). Nonlinear time series: A selective
review. \textit{Scand. J. Statist.} \textbf{21} 97--130.

\bibitem[\protect\citeauthoryear{Vasicek}{1977}]{V77}
\textsc{Vasicek, O. A.} (1977). An equilibrium characterization of the term
structure. \textit{J. Financial Economics} \textbf{5} 177--188.


\bibitem[\protect\citeauthoryear{Wu}{2003}]{wu03}
\textsc{Wu, W. B.} (2003). Empirical processes of long-memory sequences.
\textit{Bernoulli} \textbf{9} 809--831.
\MR{2047687}

\bibitem[\protect\citeauthoryear{Wu}{2005}]{wu05}
\textsc{Wu, W. B.} (2005). Nonlinear system theory: Another look at
dependence. \textit{Proc. Natl. Acad. Sci. USA} \textbf{102}
14150--14154.
\MR{2172215}

\bibitem[\protect\citeauthoryear{Wu and Mielniczuk}{2002}]{WM02}
\textsc{Wu, W. B.} and \textsc{Mielniczuk, J.} (2002). Kernel density
estimation for
linear processes. \textit{Ann. Statist.} \textbf{30} 1441--1459.
\MR{1936325}

\bibitem[\protect\citeauthoryear{Xia}{1998}]{xia98}
\textsc{Xia, Y.} (1998). Bias-corrected confidence bands in nonparametric
regression. \textit{J. R. Stat. Soc. Ser. B Stat. Methodol.}
\textbf{60} 797--811.
\MR{1649488}

\bibitem[\protect\citeauthoryear{Za\u{i}tsev}{1987}]{zait87}
\textsc{Za\u{i}tsev, A. Y.} (1987).
On the Gaussian approximation of convolutions under
multidimensional analogues of S. N. Bernstein's inequality
conditions. \textit{Probab. Theory Related Fields} \textbf{74}
535--566.
\MR{0876255}

\bibitem[\protect\citeauthoryear{Zhao}{2008}]{zhao08}
\textsc{Zhao, Z.} (2008). Parametric and nonparametric models and methods
in financial econometrics. \textit{Stat. Surv.} \textbf{2}
1--42.
\MR{2520979}

\bibitem[\protect\citeauthoryear{Zhao and Wu}{2008}]{ZW08}
\textsc{Zhao, Z.} and \textsc{Wu, W. B.} (2008).
Confidence bands in nonparametric time series regression.
\textit{Ann. Statist.} \textbf{36} 1854--1878.
\MR{2435458}

\end{thebibliography}
\end{document}